\documentclass[11pt]{article}

\usepackage[french,english]{babel}
\usepackage[T1]{fontenc}
\def\rr{{\boldsymbol{\rho}}}

\usepackage{amsfonts,amsthm}
\usepackage{amsmath,mathtools}   
\usepackage{latexsym}
\usepackage{amssymb}
\usepackage{mathrsfs}
\usepackage{pifont}
\usepackage{tikz,ulem}
\usepackage[all]{xy}
\usepackage{hyperref}
\usepackage{tikz}
\usetikzlibrary{matrix, arrows.meta, positioning}

\numberwithin{equation}{section}

\newtheoremstyle{newrem}{3pt}{3pt}{}{}% <Space above><Space below><Body font> <Indent amount>
{\bfseries}{.}{.5em}{}% <Theorem head font><Punctuation after theorem head><Space after theorem headi><Theorem head spec>

\newtheorem{theo}{Theorem}[section]

\newtheorem{rem}{Remark}[section]
\newtheorem*{theo*}{Theorem}
\newtheorem{lemm}[theo]{Lemma}
\newtheorem{prop}[theo]{Proposition}

\newtheorem{coro}[theo]{Corollary}

\theoremstyle{newrem}
\newtheorem{remax}[theo]{Remark}
\newenvironment{rema}
  {\pushQED{\qed}\remax}
  {\popQED\endremax}
\theoremstyle{definition}
\newtheorem{defi}[theo]{Definition}

\newtheorem*{term*}{Notation/Terminology}

\def \oV{{\widetilde{V}}}
\def \ow{{\widetilde{w}}}
\def \oX{{\widetilde{X}}}
\def \oZ{{\widetilde{Z}}}

\def \X{{\mathcal{X}}}
\def \Z{{\mathcal{Z}}}

\def\un{{\mathbb{I}}}
\def\V{{\mathbb{V}}}

\def\ov{{\widetilde{v}}}
\def\oalpha{{\widetilde{\alpha}}}
\def\obeta{{\widetilde{\beta}}}
\def\ogamma{{\widetilde{\gamma}}}
\def\odelta{{\widetilde{\delta}}}

\def\olambda{{\widetilde{\lambda}}}

\usepackage{geometry}
\begin{document}
\title{\bf Bispectral rational functions and Leonard trios}
\author{
Nicolas Cramp\'e\textsuperscript{$1$}\footnote{E-mail: crampe1977@gmail.com}~,
Wolter Groenevelt\textsuperscript{$2$}\footnote{E-mail:w.g.m.groenevelt@tudelft.nl}~,
Quentin Labriet\textsuperscript{$3$}\footnote{E-mail: quentin.labriet@umontreal.ca}~,
Lucia Morey\textsuperscript{$3$}\footnote{E-mail: lucia.morey@umontreal.ca}~,\\
Luc Vinet\textsuperscript{$3,4$}\footnote{E-mail: luc.vinet@umontreal.ca}~,
Carel Wagenaar\textsuperscript{$2$}\footnote{E-mail: C.C.M.L.Wagenaar@tudelft.nl}~
\vspace{0.5cm}\\
\textsuperscript{$1$}
\small CNRS -- Universit\'e de Montr\'eal CRM - CNRS,
 France.\vspace{0.2cm}\\
 \textsuperscript{$2$}
\small Technische Universiteit Delft, DIAM, PO Box 5031, 2600 GA Delft, the Netherlands.\vspace{0.2cm}\\
\textsuperscript{$3$}
\small Centre de Recherches Math\'ematiques, Universit\'e de Montr\'eal, P.O. Box 6128, \\
\small Centre-ville Station, Montr\'eal (Qu\'ebec), H3C 3J7, Canada.\vspace{0.2cm}\\
\textsuperscript{$4$}
\small IVADO, Montr\'eal (Qu\'ebec), H2S 3H1, Canada.\\
}
\date{}
\maketitle

\bigskip

\begin{center}
\begin{minipage}{12cm}
\begin{center}
{\bf Abstract}\\
\end{center}  
It is well-known that Leonard pairs have a close connection with bispectral orthogonal polynomials of the Askey scheme. In this paper, we introduce the notion of a Leonard trio $(V,\oV,Z)$, an algebraic structure extending Leonard pairs, for which the overlap coefficients of eigenfunctions of $V$ and $\oV$ are biorthogonal rational functions satisfying generalized eigenvalue problems. We introduce and start the classification of irreducible Leonard trios by using its connection with Leonard pairs and Heun operators. In particular, we show that Wilson's rational functions appear as overlap coefficients, prove its difference, recurrence and biorthogonality relations, and obtain a summation formula expressing them as a finite sum of products of two $q$-Racah polynomials. We also begin to investigate reduced Leonard trios, for which the general eigenvalue problem simplifies to a $R_I$-type recurrence relation. As an illustration, we present an example of this in which the rational functions appearing as overlap coefficients can be expressed as a ${}_{4}\phi_3$ and are associated with a Leonard pair of dual $q$-Hahn type.
\end{minipage}
\end{center}

\medskip

\section{Introduction}

Finite families of hypergeometric bispectral orthogonal polynomials and their $q$-analogs \cite{AW,Chiara,GR,Koekoek} form an important class of special
functions which is organized in part of the ($q$-)Askey
scheme. They are called bispectral since they satisfy a recurrence relation and a three-term ($q$-)difference
equation. Their classification, done in \cite{Leo}, has been enlightened  by their interpretation as overlap coefficients between two bases of the Askey--Wilson algebra or some of its limits or specialization \cite{crampe2021askey,zhedanov1991hidden} and
with the notion of the Leonard pairs introduced in \cite{Terw01,terwilliger2003introduction,TerAW,Vidunas}.
The objective of this paper consists in providing such settings for bispectral rational functions. 

A generalization of the polynomials of the ($q$-)Askey scheme to biorthogonal rational functions in terms of  balanced, very-well-poised $q$-hypergeometric functions ${}_{10}\phi_9$ has been introduced in \cite{wilson1991orthogonal,rahman1991biorthogonality,rahman1993classical,spiridonov2000spectral}. They are known as the Wilson rational functions and are bispectral satisfying a pair of generalized eigenvalue problems, as shown in \cite{zhedanov1999biorthogonal}. The algebraic interpretation of these rational functions has since been initiated for those of type Hahn \cite{tsujimoto2021algebraic} and $q$-Hahn \cite{bussiere2022bispectrality,vinet2021unified}, and even for the Wilson rational functions \cite{CTVZ}, leading to the introduction of the meta ($q$-)Hahn algebra or meta Wilson algebra. The finite-dimensional representation theory of the meta ($q$-)Hahn algebra has been carried out in \cite{tsujimoto2024meta,bernard2024meta}, providing an algebraic proof of the relations satisfied by the rational functions. Another algebraic approach to rational functions has been developed in \cite{GW2025}, which shows that (multivariate) rational functions of $q$-Racah type appear as overlap coefficients in $\mathcal{U}_q(\mathfrak{sl}_2)$.

In this paper, we introduce the definition of a Leonard trio (see Definition \ref{def:LP}), which involves three endomorphisms and two different bases. We show that the overlap coefficients between the two bases appearing in the definition satisfy two generalized eigenvalue problems, which demonstrates that these overlap coefficients are bispectral rational functions. We further prove that a biorthogonal partner also exists. The classification of Leonard trios would provide a classification of bispectral rational functions. Unfortunately, this problem seems out of reach. To focus on a restricted and tractable question, we propose the definition of irreducible Leonard trio (see Definition \ref{def:iLT}). We shall show that these more specialized trios are closely related to Leonard pairs as well as to their associated algebraic Heun operators. With this definition, we obtain the classification of the irreducible Leonard trios when the two LPs are of $q$-Racah type and prove that it is associated to the Wilson rational functions. 
We study different limits from this solution, for which another type of trios, called reduced Leonard trio (see Definition \ref{def:rLT}), appears. We show that the overlap coefficients associated with this type of trios satisfy a $R_I$-type recurrence relation (see \cite{IsMa1995}) instead of a generalized eigenvalue problem.
These results show that the Leonard trios, defined in this paper, provide a powerful algebraic framework to study bispectral rational functions. 

As a by-product of our construction, we obtain an expression of the Wilson rational function in terms of a sum of products of two $q$-Racah polynomials. This formula reduces to the Racah relation for $q$-Racah polynomials in the limit when the Wilson rational functions become the $q$-Racah polynomials.

The outline of this paper is as follows.
Section \ref{sec:LT} recalls the definition of a Leonard pair and reviews the properties relevant to this work. We then define the Leonard trio and show how it naturally leads to the consideration of biorthogonal rational functions satisfying generalized eigenvalue problems.
The definition of an irreducible Leonard trio is provided in Section \ref{sec:iLT} and its connections to Leonard pairs and the algebraic Heun operator are explained. We demonstrate that the rational functions associated with an irreducible Leonard trio can be expressed as a sum of products of two polynomials from the ($q$-)Askey scheme.
Section \ref{sec:rLT} introduces the reduced Leonard trio and establishes the constraints on the eigenvalues of its constituent operators. We show that the functions associated with a reduced Leonard trio satisfy an $R_I$-type recurrence relation rather than a generalized eigenvalue problem. Furthermore, we demonstrate that these reduced trios can be obtained as limits of irreducible ones.
Section \ref{sec:iLTqR} is dedicated to the classification of irreducible Leonard trios of $q$-Racah type. Starting from two Leonard pairs of $q$-Racah type, we derive the parameter constraints necessary to construct a Leonard trio. By solving these constraints, we show that the most general irreducible Leonard trio is associated with the Wilson rational function. As a consequence of this construction, we rederive their generalized eigenvalue problems and obtain an expression for these rational functions in terms of $q$-Racah polynomials.
Various limits of the preceding results are explored in Section \ref{sec:wrs}. We provide an example of a reduced Leonard trio and connect our findings with previously studied functions.
Finally, Section \ref{sec:out} contains concluding remarks and outlooks. Appendices \ref{sec:qRacah} and \ref{sec:Wilson} contain technical results regarding $q$-Racah polynomials and Wilson rational functions, respectively.

\section{Leonard trio \label{sec:LT}}

In this paper, $\V$ denotes a  vector space of finite dimension $N+1$ over the complex field $\mathbb{C}$, and $\mathrm{End}(\V)$ denotes the vector space over $\mathbb{C}$  of the endomorphisms of $\V$. Let $\un$ be the identity operator in $\mathrm{End}(\V)$.
A square matrix is said to be tridiagonal whenever
each nonzero entry lies on either the diagonal, the subdiagonal, or the superdiagonal. A
tridiagonal matrix is said to be irreducible whenever each entry on the subdiagonal is
nonzero and each entry on the superdiagonal is nonzero.

\subsection{Definitions \label{sec:def}}

Let us recall the definition of a Leonard pair \cite{Terw01,terwilliger2003introduction}.
\begin{defi}\label{def:LP}
 A Leonard pair (LP) is an ordered pair $(X,Y)$ of elements of $\mathrm{End(\V)}$ that satisfy the following conditions:
 \begin{itemize}
     \item[(i)]   There exists a basis of $\V$ with respect to which the matrix representing $X$ is irreducible tridiagonal, and the matrix representing $Y$ is diagonal.
     \item[(ii)]  There exists a basis of $\V$ with respect to which the matrix representing $X$ is diagonal, and the matrix representing $Y$ is irreducible tridiagonal.
 \end{itemize}
 \end{defi}
 In what follows, the notion of algebraic Heun operator will be needed.
An algebraic Heun operator associated to a LP $(X,Y)$ is the linear combination
\begin{align}\label{eq:Heun}
    H=h_0\un +h_1 X +h_2 Y +h_3 XY + h_4 YX\,,
\end{align}
 for $h_0,h_1,h_2,h_3,h_4\in \mathbb{C}$. It has been demonstrated in \cite{NT} that any endomorphism which is tridiagonal in both bases associated to the LP $(X,Y)$ can be written as \eqref{eq:Heun}. The name appears for the first time in \cite{GWZ17}, in the study of the band and time limiting problem, as a generalization of the differential Heun operator. It has also been used to simplify the computation of the entanglement entropy \cite{CNV19} and appears in the study of quantum integrable systems \cite{BP19}.

 The notion of Leonard pair has been introduced to provide  an algebraic interpretation of a theorem concerning the classification of the finite sequence of bispectral polynomials \cite{BI,Leo}.  
In what follows, we give the definition of a Leonard trio, which generalizes the notion of Leonard pair and allows us to provide a new algebraic framework for the study rational functions. We say that a diagonal matrix is multiplicity free if all entries on the diagonal are different. 
\begin{defi}\label{def:LT}
 A Leonard trio (LT) is an ordered triplet $(V,\oV,Z)$ of elements of $\mathrm{End}(\V)$ that satisfies the following properties:
   \begin{itemize}
       \item[(i)] There exists a basis of $\V$ with respect to which the matrix representing $V$ is diagonal and multiplicity free, and the matrices representing  $Z$ and $\oV Z$ are  tridiagonal.
       \item[(ii)] There exists a basis of $\V$ with respect to which the matrix representing $\oV$ is diagonal and multiplicity free, and the matrices representing  $Z$ and $Z V$ are  tridiagonal.
   \end{itemize}
\end{defi}
Let us emphasize that the matrix representing $\oV$ may be a full matrix in the first basis defined in (i) of Definition \ref{def:LT}, and $V$ may be a full matrix in the second basis defined in (ii). This is the main difference with the notion of Leonard triples introduced in \cite{Curtin} and classified in \cite{GW13,HWG13,HXG15,Huang12,KHG15}.

\begin{rem} \label{eq:affine}
From the definition of a LT, we can show these different properties:
\begin{itemize}
\item  If $(V,\oV,Z)$ is a LT then $(v_1V+v_0\un,\widetilde{v}_1\oV+\widetilde{v}_0\un,z_1Z)$, for $v_1,\widetilde{v_1},z_1\in\mathbb{C}^\star$, $v_0,\widetilde{v}_0\in \mathbb{C}$, is also a LT. We want to insist that we cannot translate $Z$ by an identity matrix.
\item
If $(V,\oV)$ is a LP, then 
    $(V,\oV,\un)$ is a LT. 
    \item   If $(V,\oV,Z)$ is a LT, then $(\oV,V,Z)$ is not necessarily a LT.
 \item If $(V,Z)$ is a LP, then $(V,V,Z)$ is a LT.   
    \end{itemize}
\end{rem}

In the remainder of this paper, for a LT $(V,\oV,Z)$, we denote by $\boldsymbol{v_n}\in\V$ and $\lambda_n\in \mathbb{C}$, for $n=0,1,\dots,N$, the eigenvectors and eigenvalues of $V$:
\begin{align}
    V \boldsymbol{v_n}= \lambda_n \boldsymbol{v_n}\,.
\end{align}
The vectors $\boldsymbol{v_n}$ form a basis of $\V$ and $\lambda_n\neq \lambda_m$ for $n\neq m$. We choose an ordering of these vectors such that, for $n=0,1,\dots, N$,
\begin{align}
     Z\boldsymbol{v_n},\ \oV Z\boldsymbol{v_n} \ \in\ \mathrm{span}_\mathbb{C}(\boldsymbol{v_{n-1}},\boldsymbol{v_n},\boldsymbol{v_{n+1}})\,,
\end{align}
where, by convention, $\boldsymbol{v_{-1}}=\boldsymbol{v_{N+1}}=0$.

We denote by $\boldsymbol{\ov_x}\in\V$ and $\olambda_x\in \mathbb{C}$, for $x=0,1,\dots,N$, the eigenvectors and the eigenvalues of $\oV$:
\begin{align}
    \oV \boldsymbol{\ov_x}=  \olambda_x \boldsymbol{\ov_x}\,.
\end{align}
The vectors $\boldsymbol{\ov_x}$ form a basis of $\V$ and $\olambda_x\neq \olambda_y$ for $x\neq y$. We choose an ordering of these vectors such that 
\begin{align}
     Z\boldsymbol{\ov_x},\ Z V\boldsymbol{\ov_x},\ \in\ \mathrm{span}_\mathbb{C}(\boldsymbol{\ov_{x-1}},\boldsymbol{\ov_x},\boldsymbol{\ov_{x+1}})\,,  
\end{align}
where, by convention, $\boldsymbol{\ov_{-1}}=\boldsymbol{\ov_{N+1}}=0$. 

\begin{rem}\label{rem:reordering}
     The reordering of the vectors $\boldsymbol{v_n}\to \boldsymbol{v_{N-n}}$ preserves the shapes of the matrices $V$, $Z$ and $\oV Z$. Similarly, the reordering of the vectors $\boldsymbol{\ov_x}\to \boldsymbol{\ov_{N-x}}$ preserves the shapes of the matrices $\oV$, $Z$ and $Z V$. 
\end{rem}

\subsection{Bispectral rational functions}

Let us endow $\V$ with a left-sesquilinear form $\langle\cdot,\cdot\rangle$ and let $\boldsymbol{v^\star_n}$ and $\boldsymbol{\ov^\star_x}$ denote the dual vectors of $\boldsymbol{v_n}$ and $\boldsymbol{\ov_x}$ defined by:
\begin{align}
   \langle\boldsymbol{v^\star_m},\boldsymbol{v_n}\rangle=\delta_{m,n}\,,\qquad \langle\boldsymbol{\ov^\star_x},\boldsymbol{\ov_y}\rangle=\delta_{x,y}\,.
\end{align}
Note that for any $\boldsymbol{w}\in\V$, we have the decompositions
\begin{align}
    \begin{split}\boldsymbol{w}&=\sum_{n=0}^N \langle \boldsymbol{v^\star_n},\boldsymbol{w}\rangle \boldsymbol{v_n}=\sum_{n=0}^N \langle \boldsymbol{v_n},\boldsymbol{w}\rangle \boldsymbol{v^\star_n}\\
    &=\sum_{x=0}^N \langle \boldsymbol{\ov^\star_x},\boldsymbol{w}\rangle \boldsymbol{\ov_x}=\sum_{x=0}^N \langle \boldsymbol{\ov_x},\boldsymbol{w}\rangle \boldsymbol{\ov^\star_x}\,.\end{split} \label{eq:basisdualbasis}
\end{align}
To any endomorphism $A\in\mathrm{End}(\V)$, we associate the adjoint endomorphism $A^\dagger\in\mathrm{End}(\V)$ defined by, for any $\boldsymbol{v},\boldsymbol{w}\in \V$,
\begin{align}
    \langle\boldsymbol{v},A\boldsymbol{w}\rangle=\langle A^\dagger \boldsymbol{v} ,\boldsymbol{w}\rangle\,.
\end{align}

The main motivation for introducing the Leonard trio is given in the following proposition, where the overlaps between the two bases of the LT provide rational functions satisfying bispectral generalized eigenvalue problems.   
\begin{prop}\label{prop:GEVPg} 
Let $(V,\oV,Z)$ be a LT.
The functions defined by
\begin{align}
    w_n(x)=\langle \boldsymbol{\ov^\star_x},\boldsymbol{v_n}\rangle\,,
\end{align}
satisfy the recurrence relation
\begin{align}
 &X_{n+1,n}\, w_{n+1}(x)  +X_{n,n}\, w_{n}(x)+ X_{n-1,n}\ w_{n-1}(x)\nonumber\\
=&\olambda_x \Big( Z_{n+1,n}\, w_{n+1}(x)+ Z_{n,n}\, w_{n}(x)    + Z_{n-1,n}\, w_{n-1}(x)\Big)\,, \label{eq:GEVPWg1}
\end{align}
and the difference equation
\begin{align}&\oX_{x,x+1}\, w_{n}(x+1)  +\oX_{x,x}\, w_{n}(x)+ \oX_{x,x-1}\ w_{n}(x-1)\nonumber\\
=&\lambda_n \Big( \oZ_{x,x+1}\, w_{n}(x+1)+ \oZ_{x,x}\, w_{n}(x)    + \oZ_{x,x-1}\, w_{n}(x-1)\Big) \,,\label{eq:GEVPWg2}
\end{align}
where $X_{m,n}$ and $Z_{m,n}$ are the entries of the matrices $\oV Z$ and $Z$ in the basis $\boldsymbol{v_n}$ and $\oX_{x,y}$ and  $\oZ_{x,y}$ are the entries of the matrices $ZV$ and $Z$ in the basis $\boldsymbol{\ov_x}$:
\begin{align}
  X_{m,n}=  \langle \boldsymbol{v_m^\star},\oV Z\boldsymbol{v_n}\rangle\,,\quad  Z_{m,n}=  \langle \boldsymbol{v_m^\star},Z\boldsymbol{v_n}\rangle\,,\quad  \oX_{x,y}=  \langle \boldsymbol{\ov_x^\star},Z V \boldsymbol{\ov_y}\rangle\,,\quad  \oZ_{x,y}=  \langle \boldsymbol{\ov_x^\star},Z  \boldsymbol{\ov_y}\rangle\,.
\end{align}
\end{prop}
\proof 
The recurrence relation \eqref{eq:GEVPWg1} is proven by computing in two different ways the following scalar product
\begin{align}
   \langle\boldsymbol{\ov^\star_x},\oV Z \boldsymbol{v_n}\rangle=\langle\oV ^\dagger \boldsymbol{\ov^\star_x},Z \boldsymbol{v_n}\rangle\,.
\end{align}
Indeed, using the property of a LT, one knows that $\oV Z$ is a tridiagonal matrix in the basis $\boldsymbol{v_n}$ which provides the lhs of \eqref{eq:GEVPWg1}.
Then, from \eqref{eq:basisdualbasis} it follows that $\oV ^\dagger \boldsymbol{\ov^\star_x}=\overline\olambda_x \boldsymbol{\ov^\star_x}$ (where $\overline{\cdot}$ stands for complex conjugation), and using the linearity of the sesquilinear form, one gets the rhs of \eqref{eq:GEVPWg1} which finishes the proof of this relation.

The difference equation is computed similarly starting from
$  \langle\boldsymbol{\ov^\star_x},ZV \boldsymbol{v_n}\rangle$, where one has to use \eqref{eq:basisdualbasis} to find the action of $Z^\dagger$ on $\boldsymbol{\ov^\star_x}$.
\endproof

The relations introduced in the previous proposition are called generalized eigenvalue problems. 
In what follows, we are going to demand that, for $0\leq n,x \leq N$, 
\begin{align}\label{eq:cond1}
    X_{n+1,n}-\olambda_x Z_{n+1,n}\neq 0\,.
\end{align}
Thanks to these conditions,
one can deduce, using recursively \eqref{eq:GEVPWg1}, that 
 $w_{n}(x)/w_0(x)$ is a rational function in terms of $\olambda_x$. 
 In addition, if $Z_{n,n}$ and $Z_{n+1,n}$ do not vanish for any $n$, one can prove by induction that $w_n(x)/w_0(x)$ is the ratio of a polynomial of degree $n$ over another one of degree $n$.
 If all the coefficients  $Z_{n+1,n}$ vanish in \eqref{eq:GEVPWg1}, then the functions $w_n/w_0$ become  polynomials, and the GEVP becomes a $R_I$ recurrence relation. 
 
We also restrict ourselves to the cases when, for $0\leq n,x \leq N$, 
\begin{align}\label{eq:cond2}
    \oX_{x,x+1}-\lambda_n \oZ_{x,x+1}\neq 0\,,
\end{align}
 such that relation \eqref{eq:GEVPWg2} can be also used recursively to define $w_{n}(x+1)$ in terms of $w_n(0)$. 
In other words, the previous proposition, together with the two constraints \eqref{eq:cond1} and \eqref{eq:cond2}, means that the overlap coefficients between the two bases introduced in the definition of a LT are bispectral rational functions. It generalizes the fact that the overlap coefficients between the two bases introduced in the definition of a LP are bispectral polynomials.\\  
   
The rational functions $w_n$, introduced in the previous proposition, have natural partners to obtain a biorthogonality relation given in the next proposition.
\begin{prop}\label{pro:bpw}
Define $w_n(x)$ as in Proposition \ref{prop:GEVPg} and the functions $\ow_n(x)$ as
\begin{align}
    \ow_n(x)=\langle \boldsymbol{v^\star_n},\boldsymbol{\ov_x}\rangle\,.
\end{align}
These functions are biorthogonal and close to identity
\begin{align}\label{eq:biort}
    \sum_{x=0}^N w_m(x)\ow_n(x)=\delta_{m,n}\,,\qquad
    \sum_{n=0}^N w_n(x)\ow_n(y)=\delta_{x,y}\,.
\end{align}
\end{prop}
\proof
Using \eqref{eq:basisdualbasis}, one can express $\boldsymbol{v_m}$ in the basis $\boldsymbol{\ov_x}$ to obtain 
\[
    \delta_{m,n}=\langle \boldsymbol{v^\star_n} ,\boldsymbol{v_m}\rangle= \sum_{x=0}^N\langle \boldsymbol{\ov^\star_x} ,\boldsymbol{v_m}\rangle\langle \boldsymbol{v^\star_n} ,\boldsymbol{\ov_x}\rangle  \,, 
\]
proving the first relation in \eqref{eq:biort}. The second relation is proven similarly.
\endproof
At this point, the biorthogonal partners $\ow_n(x)$ do not necessarily satisfy recurrence or difference relations, although this will be the case for the examples we will discuss later.

 The next objective is to obtain a classification of LTs in analogy with the classification of LPs \cite{Terw01,terwilliger2003introduction}. However, this question seems out-of-reach at this moment. By imposing a few additional conditions, one arrives at the definition of an irreducible Leonard trio, where we succeed to compute a large class of solutions in the following sections.

\section{Irreducible Leonard trio \label{sec:iLT}}

\subsection{Definition and algebraic  Heun operators}

To simplify the classification but still have large classes of examples, we introduce the notion of irreducible Leonard trio.
\begin{defi}\label{def:iLT}
 An irreducible Leonard trio (iLT) is an ordered triplet $(V,\oV,Z)$ of elements of $\mathrm{End}(\V)$ satisfying the following properties:
   \begin{itemize}
       \item[(i)] There exists a basis of $\V$ with respect to which the matrix representing $V$ is diagonal, and the matrix representing $\oV Z$ is  tridiagonal and the one representing $Z$ is irreducible tridiagonal.
       \item[(ii)] There exists a basis of $\V$ with respect to which the matrix representing $\oV$ is diagonal and the matrix representing $ Z V$ is  tridiagonal and the one representing $Z$ is irreducible tridiagonal.
       \item[(iii)] There exists a basis of $\V$ with respect to which the matrix representing $Z$ is diagonal and the matrices  representing $V$ and $\oV$ are irreducible tridiagonal.
   \end{itemize}
\end{defi}
For this definition, there exist strong connections with Leonard pairs as stated in the following proposition.
\begin{prop}
    Let $(V,\oV,Z)$ be an iLT. The following statements hold:
    \begin{itemize}
        \item $(V,Z)$ is a LP;
        \item $(\oV,Z)$ is a LP;
        \item $(V,\oV,Z)$ is a LT.
    \end{itemize}    
\end{prop}
\proof The proof is a direct consequence of the different definitions except for the non-degeneracy of the spectrum of $V$ and $\oV$ in a LT. This follows from the fact that the two elements of a LP have a non-degenerate spectrum \cite[Lemma 3.1]{TerAW}. \endproof

For an iLT, the matrices representing $Z$ in the two bases associated to the LT are irreducible tridiagonal. It explains the choice of the adjective \textit{irreducible} to characterize this type of LT. The iLT has the following interesting connection with algebraic Heun operators (see below Definition \ref{def:LP} for their definition) which is very helpful for their classification.  
\begin{prop}\label{pro:Heun}
    Let  $(V,\oV,Z)$ be an iLT. There exist ${h}_0,h_1,\dots,{h}_9\in\mathbb{C}$ such that
   \begin{align}\label{eq:Heung1}
    \oV Z={h}_0\un+{h}_1Z+{h}_2 V+{h}_3VZ+{h}_4ZV\,,\\
    \label{eq:Heung4}
    Z V =h_5 \un+h_6Z+h_7 \oV+h_8\oV Z+h_9Z\oV\,.
\end{align}
\end{prop}
\proof
By the definition of a LP, the element $\oV$ acts tridiagonally on a basis of $\V$ with respect to which $Z$ is diagonal. 
Then, the element $\oV Z$ acts tridiagonally on this basis.
In addition, the operator $\oV Z$ is also tridiagonal in the basis $\boldsymbol{v_n}$ (by definition of a LT). Hence from \cite{NT} we know that $\oV Z$, which is tridiagonal in the eigenbases of the two operators $V$ and $Z$ defining the LP $(V,Z)$, must be a Heun operator associated to that Leonard pair \textit{i.e.} of the form \eqref{eq:Heung1}. 
Similarly, we can prove that $Z V$ is an algebraic Heun operator with respect to $\oV$ and $Z$. \endproof

A reverse of the Proposition \ref{pro:Heun} is also true.
\begin{prop}\label{pro:reverse}
    Let  $(V,Z)$ and $(\oV,Z)$ be two LPs. If there exist ${h}_0,h_1,\dots,{h}_9\in\mathbb{C}$ such that
   \eqref{eq:Heung1} and \eqref{eq:Heung4} hold, then $(V,\oV,Z)$ is an iLT.
\end{prop}
\proof The points (i) and (ii) of the definition of a LT is obtained using the expressions of $\oV Z$ and $ Z V$ in terms of the algebraic Heun  operators. The other points are directly obtained by the hypothesis of the proposition. \endproof

Using Proposition \ref{pro:Heun}, we can show that the biorthogonal partner $\ow_n(x)$ satisfies recurrence and difference relations as well. Indeed, from \eqref{eq:Heung1}, we get that
\begin{align*}
    h_4 ZV=\oV Z-{h}_0\un-{h}_1Z-{h}_2 V-{h}_3VZ\,.
\end{align*}
Multiplying the lhs of \eqref{eq:Heung4} by $h_4$ and substituting above expression gives,
\begin{align*}
    \oV Z-{h}_0\un-{h}_1Z-{h}_2 V-{h}_3VZ =h_4h_5 \un+h_4h_6Z+h_4h_7 \oV+h_4h_8\oV Z+h_4h_9Z\oV\,.
\end{align*}
Rearranging terms gives that, besides $ZV$, also $V(Z+h_2/h_3 \un)$ is a Heun operator associated to $(Z,\oV)$. Similarly, $(Z+h_7/h_9 \un)\oV$ is a Heun operator associated to $(Z,V)$.
\begin{coro} \label{cor:otherHeun}
    Let ${h}_0,h_1,\dots,{h}_9\in\mathbb{C}$ as in Proposition \ref{pro:Heun} and assume $h_3$ and $h_9$ are nonzero, then
    \begin{align}
        \big(Z+\frac{h_7}{h_9} \un\big)\oV &= -\frac{h_5+h_0h_8}{h_9}\un -\frac{(h_6+h_1h_8)}{h_9} Z -\frac{h_2h_8}{h_9} V -\frac{h_3h_8}{h_9} VZ + \frac{(1-h_4h_8)}{h_9} ZV \,, \\
        V\big(Z+\frac{h_2}{h_3} \un\big) &=-\frac{({h}_0+h_4h_5)}{h_3}\un-\frac{({h}_1+h_4h_6)}{h_3}Z -\frac{h_4h_7}{h_3} \oV+\frac{(1-h_4h_8)}{h_3}\oV Z-\frac{h_4h_9}{h_3}Z\oV\,.
    \end{align}
\end{coro}
Using the above corollary, one can deduce recurrence and difference relations for $\ow_n(x)$ using the same method as in Proposition \ref{prop:GEVPg}.\\

In addition to the connection with algebraic Heun operators, the fact that iLTs are related to Leonard pairs allows us to use the classification of the latter to explore the possible types of iLTs. It is known that any Leonard pair is associated with a family of bispectral polynomials. These families must belong to the finite sequences of the $q$-Askey scheme ($q$-Racah, $q$-Hahn, dual $q$-Hahn, $q$-Krawtchouk, dual $q$-Krawtchouk, affine $q$-Krawtchouk, quantum $q$-Krawtchouk), the limit $q\to 1$ of that scheme (Racah, Hahn, dual Hahn, Krawtchouk) or its $q\to -1$ limit (Bannai--Ito and their descendants).
In the following subsection, we recall the precise relation between LPs and bispectral polynomials.

\subsection{Leonard pairs and bispectral polynomials \label{sec:LPbp}}

If the finite sequence of polynomials $(P_n)_{n=0}^N$ is one of the sequences mentioned above, then they satisfy a recurrence relation, for $0\leq n,x\leq N$,
     \begin{align}
     \lambda_x P_n(x)&=\mathcal{A}_n P_{n+1}(x)+\mathcal{B}_nP_n(x)+\mathcal{C}_nP_{n-1}(x)\,.\label{eq:recurrence-P}
     \end{align}
By the Favard theorem \cite{Chiara}, these polynomials are orthogonal, \textit{i.e.} there exists weight $\Omega_{x}$ such that the following relations hold, for $0\leq m,n\leq N$,
\begin{align}
    \sum_{x=0}^N \Omega_{x}P_m(x)P_n(x)=\omega_n \delta_{m,n}\,.
\end{align}
In addition to these relations, they also satisfy difference equations, for $0\leq n,x\leq N$,
      \begin{align}
\xi_n P_n(x)&=\mathscr{A}_xP_n(x+1)+\mathscr{B}_xP_n(x)+\mathscr{C}_xP_n(x-1)\,,\label{eq:difference-P}
     \end{align}
and close to identity
\begin{align}\label{eq:cig}
     \sum_{n=0}^N \frac{1}{\omega_n}P_n(x)P_n(y)=\frac{1}{\Omega_x}\delta_{x,y}\,.
\end{align}
To connect these polynomials to LPs, let $\boldsymbol{z_i}$ ($i=0,1,\dots,N$) be a basis of $\V$ and equip $\V$ with the inner product induced by $\langle \boldsymbol{z_i},\boldsymbol{z_j}\rangle=\delta_{i,j}$. Define two operators $Z$ and $V$ by, for $i=0,1,\dots,N$,
\begin{align}
Z\boldsymbol{z_i}=\xi_i \boldsymbol{z_i}\,, 
\qquad
    V\boldsymbol{z_i}&=\mathcal{C}_{i+1}\boldsymbol{z_{i+1}}+\mathcal{B}_i\boldsymbol{z_i}
   +\mathcal{A}_{i-1}\boldsymbol{z_{i-1}}\,.
   \end{align}
These two operators form a LP $(V,Z)$. Indeed, defining $\boldsymbol{v_n}$ by, for $n=0,1,\dots,N$,
\begin{align}\label{eq:vnzi-general}
    \boldsymbol{v_n}=\sum_{i=0}^N P_i(n)\ \boldsymbol{z_{i}}\,,
\end{align}
one can show, using the recurrence relation \eqref{eq:recurrence-P}, that $V$ acts diagonally on this basis,
\begin{align}
    V \boldsymbol{v_n}=\lambda_n\boldsymbol{v_n}\,,
\end{align}
and, using the difference relation  \eqref{eq:difference-P}, that $Z$ acts tridiagonally:
\begin{align}
  &   Z \boldsymbol{v_n} =\mathscr{A}_{n} \boldsymbol{v_{n+1}}+\mathscr{B}_{n}\boldsymbol{v_n}+\mathscr{C}_{n}\boldsymbol{v_{n-1}}\,.
\end{align}
Using relation \eqref{eq:cig}, the dual vectors $\boldsymbol{v^\star_n}$ can be written as follows:
\begin{align}
    \boldsymbol{v^\star_n}=\sum_{i=0}^N \frac{\omega_i}{\Omega_n} P_i(n)\boldsymbol{z_i}\,.
\end{align}
We deduce that the overlap coefficients between the two bases $\boldsymbol{z_i}$, $\boldsymbol{v_n}$ are proportional to the polynomials
\begin{align}
 \langle \boldsymbol{z_i},\boldsymbol{v_n} \rangle= P_i(n)\,,\qquad   \langle \boldsymbol{v^\star_n},\boldsymbol{z_i} \rangle=\frac{\omega_i}{\Omega_n} P_i(n)\,.
\end{align}

\subsection{Irreducible Leonard trios and bispectral polynomials \label{sec:LTbp}}

Let  $(V,\oV,Z)$ be an irreducible Leonard trio. By definition, $(V,Z)$ and $(\oV,Z)$ are Leonard pairs, and we can associate to them (as recalled in the previous subsection)  the finite sequences of polynomials $(P_i)_{i=0}^N$ and $(\widetilde{P}_i)_{i=0}^N$, respectively.
For the polynomials $P_n$, we use the notations of the previous subsection whereas for $\widetilde P_n$, we use the same names for all the functions appearing in the relations satisfied by these polynomials but with an additional tilde. For example, the vector of the basis diagonalizing $\oV$ is denoted by $\boldsymbol{\ov_n}$ and the equivalent of relation \eqref{eq:vnzi-general} reads
\begin{align}\label{eq:ovxvi7-general}
    \boldsymbol{\ov_x}=\sum_{i=0}^N \widetilde{P}_i(x)\boldsymbol{z_{i}}\,.
\end{align}

From these results, we can obtain a relation between the rational function $w_n(x)$, defined in Proposition \ref{prop:GEVPg}, and both sets of polynomials associated to both LPs,
\begin{align}\label{eq:wWilsong}
w_{n}(x)=\langle \boldsymbol{\ov_x^\star},\boldsymbol{v_n}\rangle=\sum_{i=0}^N\langle \boldsymbol{\ov_x^\star},\boldsymbol{z_i}\rangle\langle \boldsymbol{z_i},\boldsymbol{v_n}\rangle= \sum_{i=0}^N \frac{\widetilde{\omega}_x}{\widetilde{\Omega}_i}\widetilde{P}_i(x)P_i(n)\,.
\end{align}
We remark that this construction shows that the bispectral rational functions $w_n(x)$ we are interested in are sums of the product of two bispectral polynomials.

As stated in Proposition \ref{pro:Heun}, associated to an iLT there are the algebraic Heun operators \eqref{eq:Heung1} and \eqref{eq:Heung4}. Using that in the basis $\boldsymbol{z_i}$, the operators $V$ and $\oV$ are tridiagonal and $Z$ is diagonal, we can compute the constants $h_i$ appearing in these relations.
In addition, if we know the expressions of the parameters $h_i$, we can determine the action of $\oV Z$ on $\boldsymbol{v_n}$ using  \eqref{eq:Heung1}
\begin{align}
  (\oV Z)\boldsymbol{v_n}=&(h_0+h_2\lambda_n +(h_1+h_3\lambda_n+h_4\lambda_n)\mathscr{B}_n) \boldsymbol{v_n}\nonumber\\
 \label{eq:www1}
+&(h_1+h_3\lambda_{n+1}+h_4\lambda_n)\mathscr{A}_n\boldsymbol{v_{n+1}}
+(h_1+h_3\lambda_{n-1}+h_4\lambda_n)\mathscr{C}_n\boldsymbol{v_{n-1}}\,,
\end{align}
and we can compute the action of $ZV$ on $\boldsymbol{\ov_x}$ using  \eqref{eq:Heung4}
\begin{align}
  ( ZV)\boldsymbol{\ov_x}=&(h_5+h_7\olambda_x +(h_6+h_8\olambda_x+h_9\olambda_x)\widetilde{\mathscr{B}}_x) \boldsymbol{\ov_x}\nonumber\\
   \label{eq:www2}
+&(h_6+h_8\olambda_{x+1}+h_9\olambda_x)\widetilde{\mathscr{A}}_x\boldsymbol{\ov_{x+1}}
+(h_6+h_8\olambda_{x-1}+h_9\olambda_x)\widetilde{\mathscr{C}}_x\boldsymbol{\ov_{x-1}}\,.
\end{align}
 These formulas allow us to obtain the explicit expressions of the functions appearing in the GEVP \eqref{eq:GEVPWg1} and \eqref{eq:GEVPWg2} for $w_n(x)$:
 \begin{subequations}
  \begin{align}
  &  Z_{n+1,n}=\mathscr{A}_n\,,\quad Z_{n,n}=\mathscr{B}_n\,,\quad Z_{n-1,n}=\mathscr{C}_n\,,\\
  &  X_{n+1,n}=(h_1+h_3\lambda_{n+1}+h_4\lambda_n)\mathscr{A}_n\,,\quad X_{n,n}=h_0+h_2\lambda_n +(h_1+h_3\lambda_n+h_4\lambda_n)\mathscr{B}_n\,,\\
  &X_{n-1,n}=(h_1+h_3\lambda_{n-1}+h_4\lambda_n)\mathscr{C}_n\,,
\end{align}   
 \end{subequations}
and
\begin{subequations}
    \begin{align}
  &  \oZ_{x,x-1}=\widetilde{\mathscr{A}}_{x-1}\,,\quad \oZ_{x,x}=\widetilde{\mathscr{B}}_x\,,\quad \oZ_{x,x+1}=\widetilde{\mathscr{C}}_{x+1}\,,\\
  &  \oX_{x,x-1}=(h_6+h_8\olambda_{x-1}+h_9\olambda_x)\widetilde{\mathscr{A}}_{x-1}\,,\quad X_{x,x}=h_5+h_7\olambda_x +(h_6+h_8\olambda_x+h_9\olambda_x)\widetilde{\mathscr{B}}_x\,,\\
  &\oX_{x,x+1}=(h_6+h_8\olambda_{x}+h_9\olambda_{x+1})\widetilde{\mathscr{C}}_{x+1}\,.
\end{align}
\end{subequations}

Conversely, to classify the iLT, one can start with two LPs $(V,Z)$ and $(\oV,Z)$ (or the two associated sequences of polynomials) and try to find out when the triplet $(V,\oV,Z)$ forms an iLT.
Let us emphasize that when one considers only one LP $(V,Z)$, any affine transformations of $V$ and $Z$ or normalizations of the vectors of the two bases associated to the LP, only change the bispectral polynomials in a trivial way. In the case of an iLT, the two LPs share one element: the operator $Z$. Therefore, we can choose, for example, the normalizations for the vectors $\boldsymbol{z_{i}}$ such that $V$ takes the usual form, but the action of $\oV$ may involve a normalization resulting in an extra parameter of the rational function $w_n(x)$. That is, the term $\widetilde{P_i}(x)$ in \eqref{eq:ovxvi7-general} may contain a factor depending on $i$ and a new parameter, having a non-trivial influence on the rational function $w_n(x)$ in \eqref{eq:wWilsong}. 
In addition, as stated in Remark \ref{eq:affine}, 
 for $(V,\oV,Z)$ a LT,  $(V,\oV,Z+z_0\un)$, ($z_0\in \mathbb{C}$), is not necessarily a LT. This implies that we cannot choose freely this constant when we define $Z$ (see for example the constant $\sigma$ in relation \eqref{eq:Zzi}).
When these two LPs are chosen, we can express the relations \eqref{eq:Heung1} and \eqref{eq:Heung4} in the basis $\boldsymbol{z_i}$. This gives constraints on the parameters $h_i$ as well as on the parameters of the polynomials. Finding all the solutions of these constraints provides a classification of the iLT associated with the two polynomials we start with.
In section \ref{sec:iLTqR}, we follow this program to classify the iLT when $(V,Z)$ and $(\oV,Z)$ are both associated to the most general finite bispectral polynomials connected to LPs, namely the $q$-Racah polynomials. 

\section{Reduced Leonard trio \label{sec:rLT}}

We have just defined the irreducible Leonard trio but 
we can also consider a slightly different simplification of a LT. In the following, an irreducible lower (resp. upper) bidiagonal matrix is a lower (resp. upper) bidiagonal matrix with the non-vanishing entries on the lower (resp. upper) diagonal.
\begin{defi}\label{def:rLT}
 An upper (resp. lower) reduced Leonard trio (rLT) is an ordered triplet $(V,\oV,Z)$ of elements of $\mathrm{End}(\V)$ that satisfies the following properties:
   \begin{itemize}
       \item[(i)] There exists a basis of $\V$ with respect to which the matrix representing $V$ is diagonal, the one representing $\oV Z$ is tridiagonal, and the one representing $Z$ is irreducible upper bidiagonal (resp. tridiagonal).
       \item[(ii)] There exists a basis of $\V$ with respect to which the matrix representing $\oV$ is diagonal, the one  representing $ Z V$ is  tridiagonal, and the one representing $Z$ is irreducible tridiagonal (resp. lower bidiagonal).
       \item[(iii)] There exists a basis of $\V$ with respect to which the matrix representing $Z$ is diagonal, the one representing  $\oV$ is irreducible tridiagonal (resp. lower bidiagonal), and $V$ is irreducible upper bidiagonal (resp. tridiagonal).
   \end{itemize}
\end{defi}
 A pair $(V,Z)$ which follows the definition of a LP but with the properties to be tridiagonal replaced by the fact to be lower (resp. upper) bidiagonal is called a reduced Leonard pair. 
In the following, such pairs appear as a certain limit of LP. 
\begin{rema}
If $(V,\oV,Z)$ is an upper rLT then the pair $(\oV,Z)$ is a LP and $(V,Z)$ is a reduced LP .
If $(V,\oV,Z)$ is a lower rLT then the pair $(V,Z)$ is a LP and $(\oV,Z)$ is a reduced LP.
\end{rema}

In the above definition, we made a choice of lower (or upper) bidiagonal matrices. This choice seems arbitrary but we can obtain a lower matrix from an upper one by reversing the order of the vectors in the corresponding basis. We made this choice such that the corresponding overlap coefficients have nice properties as shown below.

\subsection{Constraints for the eigenvalues of $V$ and $Z$ when $(V,Z)$ is a reduced LP \label{sec:consrLP}} 

Let $(V,\oV,Z)$ be a reduced Leonard trio with $(V,Z)$ a reduced LP.
As previously, let us denote by $\boldsymbol{z_i}$ the eigenvectors of $Z$ and $\xi_i$, the associated eigenvalues. In this basis, $V$ is upper bidiagonal with non-zero entries on the upper diagonal. Therefore, the action of $V$ can be written as
\begin{align}
    V \boldsymbol{z_i}=\lambda_i \boldsymbol{z_i} +{\mathcal{A}}_{i-1}\boldsymbol{z_{i-1}}\,,
\end{align}
and the eigenvectors $\boldsymbol{v_x}$ of $V$ are 
\begin{align}\label{eq:rLTvx}
    \boldsymbol{v_x}=\sum_{i=0}^x \zeta(x,i) \boldsymbol{z_i}\,,
\end{align}
with $\zeta$ satisfying, for $i<x$,
\begin{align}\label{eq:xi1}
    (\lambda_x-\lambda_i)\zeta(x,i)=\zeta(x,i+1){\mathcal{A}}_{i}\,.
\end{align}
The action of $Z$ is of the form
\begin{align}
    Z\boldsymbol{v_x}=\xi_x\boldsymbol{v_x}+{\mathscr{C}}_x\boldsymbol{v_{x-1}}\,,
\end{align}
and using the expression \eqref{eq:rLTvx}, the following constraint must be verified
\begin{align}\label{eq:xi2}
    (\xi_i-\xi_x)\zeta(x,i)=\zeta(x-1,i){\mathscr{C}}_x\,. 
\end{align}
Relations \eqref{eq:xi1} and \eqref{eq:xi2} allow to express recursively $\zeta(x,i)$ in terms of $\zeta(0,0)$ but there are different paths to obtain this result, and they must be compatible. These compatibility conditions are obtained by comparing the two ways to express $\zeta(x,i+1)$ in terms of $\zeta(x-1,i)$:
\begin{align}\label{eq:ltzeta}
  \zeta(x,i+1)&= \frac{\lambda_x-\lambda_i}{{\mathcal{A}}_{i}} \zeta(x,i)=\frac{{\mathscr{C}}_{x}(\lambda_x-\lambda_i)}{{\mathcal{A}}_{i}(\xi_i-\xi_x)} \zeta(x-1,i)\\
  &=\frac{{\mathscr{C}}_{x}}{\xi_{i+1}-\xi_x} \zeta(x-1,i+1)=\frac{{\mathscr{C}}_{x}(\lambda_{x-1}-\lambda_i)}{{\mathcal{A}}_{i}(\xi_{i+1}-\xi_x)} \zeta(x-1,i)\,,
\end{align}
which leads to 
\begin{align}
  \frac{\lambda_x-\lambda_i}{\lambda_{x-1}-\lambda_i} =\frac{\xi_x-\xi_i}{\xi_x-\xi_{i+1}}\,.  
\end{align}
The eigenvalues $\xi_x$ of $Z$, which is a member of the LP $(\oV,Z)$, can take only three different forms \cite{Terw01}. The first form is  $\xi_x=\mathfrak{a}_0+\mathfrak{a}_1 q^x+\mathfrak{a}_2 q^{-x}$. Using relation \eqref{eq:ltzeta} for $i=0$, the eigenvalues take necessarily the form
   $\lambda_x= \mathfrak{a}_3+\mathfrak{a}_4\frac{1-q^x}{\mathfrak{a}_1q^{x+1}-\mathfrak{a}_2}$. Injecting this form in \eqref{eq:ltzeta}, it imposes that $\mathfrak{a}_1=0$ or $\mathfrak{a}_2=0$. Therefore, in this case, only the two following sequences of eigenvalues are possible
   \begin{align}\label{eq:qr}
    (\xi_x,\lambda_x)=(\mathfrak{a}_0+\mathfrak{a}_1 q^x,\mathfrak{a}_3+\mathfrak{a}_4q^{-x})\,,\quad \text{or}\qquad (\mathfrak{a}_0+\mathfrak{a}_2 q^{-x},\mathfrak{a}_3+\mathfrak{a}_4q^{x})\,.   \end{align}
    In the previous relations, we rescale the parameters $\mathfrak{a}_3$ and $\mathfrak{a}_4$ to simplify the presentation.
    Starting from the second possibility $\xi_x=\mathfrak{a}_0+\mathfrak{a}_1 x+\mathfrak{a}_2 x^2$ similar computations lead to only one solution
   \begin{align}\label{eq:lr}
       (\xi_x,\lambda_x)= (\mathfrak{a}_0+\mathfrak{a}_1 x,  \mathfrak{a}_3+\mathfrak{a}_4 x)\,.
   \end{align}
From the third possibility $\xi_x=\mathfrak{a}_0+\mathfrak{a}_1 (-1)^x+\mathfrak{a}_2 (-1)^x x$ ($\mathfrak{a_2}\neq 0$), there does not exist any solution. 

Let us emphasize that there are no constraints on ${\mathcal{A}}_i$ and ${\mathscr{C}}_x$ to be a reduced LP. Some constraints will appear when they are a part of a rLT.
To conclude, we show that the only possibilities for the eigenvalues of $Z$ and $V$ are given by \eqref{eq:qr}  and \eqref{eq:lr}. We denote by $q$-red (resp. $l$-red) a reduced Leonard pair whose eigenvalues satisfy \eqref{eq:qr} (resp. \eqref{eq:lr}). 

Same constraints for the eigenvalues of $\oV$ and $Z$,
when $(V,\oV,Z)$ is a reduced Leonard trio with $(\oV,Z)$ a reduced LP, can be found.

\subsection{Overlap coefficients for reduced LT} 

Let $(V,\oV,Z)$ be a reduced Leonard trio  with $(V,Z)$ a reduced LP. The overlap coefficients of the LP $(\oV,Z)$ is still expressed in terms of the  bispectral orthogonal polynomials $\widetilde{P}_i$:
\begin{align}
 \langle  \boldsymbol{\ov^\star_x}, \boldsymbol{z_i} \rangle= \frac{\widetilde{\omega}_x}{\widetilde{\Omega}_i}\widetilde{P}_i(x)\,.
\end{align}
Solving \eqref{eq:xi1}, one obtains $\zeta(x,i)=\prod_{j=0}^{i-1}\frac{\lambda_x-\lambda_j}{{\mathcal{A}}_j}$, ($i\leq x$) and the overlap coefficients $w_n(x)$ associated to this rLT  is given by
\begin{align}
    w_n(x)=\langle\boldsymbol{\ov^\star_x} , \boldsymbol{v_n} \rangle
    =\sum_{i=0}^N  \frac{\widetilde{\omega}_x}{\widetilde{\Omega}_i}\widetilde{P}_i(x) \prod_{j=0}^{i-1}\frac{\lambda_x-\lambda_j}{{\mathcal{A}}_j}\,.
\end{align}
In this case, by definition, $Z$ is upper bidiagonal therefore the coefficients $Z_{n+1,n}$ vanish. 
The recurrence GEVP \eqref{eq:GEVPWg1} simplifies to
\begin{align}
 &X_{n+1,n}\, w_{n+1}(x)  +X_{n,n}\, w_{n}(x)+ X_{n-1,n}\ w_{n-1}(x)
=\olambda_x \Big( Z_{n,n}\, w_{n}(x)    + Z_{n-1,n}\, w_{n-1}(x)\Big)\,, \label{eq:RIWg1}
\end{align}
and becomes a $R_I$ recurrence relation.

For the case when $(V,\oV,Z)$ is a reduced Leonard trio  with $(\oV,Z)$ a reduced LP, it is the difference GEVP \eqref{eq:GEVPWg2} which simplifies and becomes of type $R_I$. 

\subsection{Limits for LP of $q$-Racah type and reduced LP}
In this section we show how one can obtain a reduced LP from a LP of $q$-Racah type. To do so, we consider the $q$-Racah polynomial $R_n^{(qR)}$ (see
Appendix \ref{sec:qRacah} for definition and properties) and the functions 
 $r_n(x;\alpha,\gamma,\delta)$ defined by, for $\gamma=q^{-N-1}$,
\begin{align}
  r_n(x;\alpha,\gamma,\delta) =\mathop{\text{lim}}_{\beta \to 0} \beta^n R_n^{(qR)}(x;\alpha,\beta,\gamma,\delta/\beta) \,.
\end{align}
Using \eqref{eq:lim0}, one obtains that
\begin{align}
     r_n(x;\alpha,\gamma,\delta) =
        \displaystyle  \frac{(q^{-x};q)_n}{(\alpha q,\gamma q,\delta q;q)_n }\left(\gamma \delta q^{x+1}\right)^n.
\end{align}
Let us remark that $r_n(x)$ vanishes for $n>x$. Taking this limit for the recurrence relation and the difference equation of the $q$-Racah polynomials (see \eqref{eq:recuR} and \eqref{eq:diffR}), one obtains the following relations for $r_n$: 
\begin{align}
&\gamma\delta q^{x+1}r_n(x)=(1-\alpha q^{n+1})(1-\delta q^{n+1})(1-\gamma q^{n+1})r_{n+1}(x)+\gamma\delta q^{n+1} r_n(x)\,,\\
 &  q^{-n}r_n(x)=  (1-q^{-x}) r_n(x-1) +q^{-x}r_n(x) \,.
\end{align}
Following the steps done in Subsection \ref{sec:LPbp} for the LP, one can define two operators $V$ and $Z$ by their actions on the vectors $\boldsymbol{z_i}$:
\begin{align}
    Z\boldsymbol{z_i}=q^{-i}\boldsymbol{z_i} \,,\qquad
    V\boldsymbol{z_i}=(1-\alpha q^{i})(1-\delta q^{i})(1-\gamma q^{i})\boldsymbol{z_{i-1}} +\gamma\delta q^{i+1} \boldsymbol{z_{i}} \,.
\end{align}
Defining new vectors by 
\begin{align}
    \boldsymbol{v_n}=\sum_{i=0}^N r_i(n;\alpha,\gamma,\delta)  \boldsymbol{z_{i}}\,,
\end{align}
one can show, using the previous relation, that
\begin{align}
     Z\boldsymbol{v_n}=(1-q^{-n})\boldsymbol{v_{n-1}}+q^{-n}\boldsymbol{v_n} \,,\qquad
    V\boldsymbol{v_n}=\gamma\delta q^{n+1}  \boldsymbol{v_n}\,.
\end{align}
We can conclude that $(V,Z)$ forms a reduced LP. Let us remark that the eigenvalues of this pair are compatible with the results found in Subsection \ref{sec:consrLP}.

\section{Irreducible Leonard trio of $q$-Racah type \label{sec:iLTqR}} 

This section is devoted to providing the classification of the irreducible Leonard trio of $q$-Racah type following the procedure explained at the end of Subsection \ref{sec:LTbp}.
Solutions associated with other families can be obtained by different limits, as shown in the following section.

\subsection{Constraints to obtain a LT}

Define the operators $Z$, $V$ and $\oV$ in the basis $\boldsymbol{z_i}$ in which $Z$ acts diagonally and such that $(V,Z)$, $(\oV,Z)$ are LPs of type $q$-Racah: 
\begin{align}
 Z\boldsymbol{z_i}&=\underbrace{(q^{-i}+\alpha\beta q^{i+1}+\sigma)}_{=\zeta_i}
\boldsymbol{z_i}\,,\label{eq:Zzi}\\
    V\boldsymbol{z_i}&=C^{(qR)}_{i+1}(\alpha,\beta,\gamma,\delta)\boldsymbol{z_{i+1}}+B^{(qR)}_i(\alpha,\beta,\gamma,\delta)\boldsymbol{z_i}
   +A^{(qR)}_{i-1}(\alpha,\beta,\gamma,\delta)\boldsymbol{z_{i-1}}\,,\label{eq:Vz}\\
   \oV\boldsymbol{z_i}&=C^{(qR)}_{i+1}(\oalpha,\obeta,\ogamma,\odelta)\nu_{i+1}\boldsymbol{z_{i+1}}+B^{(qR)}_i(\oalpha,\obeta,\ogamma,\odelta)\boldsymbol{z_i}+A^{(qR)}_{i-1}(\oalpha,\obeta,\ogamma,\odelta)\frac{1}{\nu_i}\boldsymbol{z_{i-1}} \,,\label{eq:Vtz}
   \end{align}
   where $A^{(qR)}$, $B^{(qR)}$ and $C^{(qR)}$ are defined in Appendix \ref{sec:qRacah}.
 We already use the different symmetry explained in Remark \ref{eq:affine} to reduce the number of possible parameters. 
 For example, we use the symmetry $V\to v_1V+v_0\un$ to choose the global factor in the action of $V$ and to get that $B_i^{(qR)}$ is given by \eqref{eq:BAC}. As explained at the end of Subsection \ref{sec:LTbp}, after these simplifications, it still remains to find the parameter $\sigma$ in the eigenvalue of $Z$ and the normalization factor $\nu_i$ in the action of $\oV$.
Since the dimension of the operators must be the same and $V$ and $\oV$ share the same partner $Z$, the parameters must satisfy:
\begin{align}
 \obeta=\alpha\beta/\oalpha\,,\qquad \ogamma=\gamma=q^{-N-1}\,.
\end{align} 
Let us mention that,
instead of considering action \eqref{eq:Vtz}, we could have considered the following action: \begin{align}
    \oV\boldsymbol{z_i}&=C^{(qR)}_{N-i+1}(\oalpha',\obeta',\ogamma',\odelta')\nu'_{i-1}\boldsymbol{z_{i-1}}+B^{(qR)}_{N-i}(\oalpha',\obeta',\ogamma',\odelta')\boldsymbol{z_i}+A^{(qR)}_{N-i-1}(\oalpha',\obeta',\ogamma',\odelta')\frac{1}{\nu'_{i}}\boldsymbol{z_{i+1}} \,.\label{eq:Vtzbis}
\end{align} 
However setting $\oalpha'=\ogamma/\obeta$, $\obeta'=\ogamma/\oalpha$, $\ogamma'=\ogamma$, $\odelta'=\odelta$ and $\nu'_i=1/\nu_{i+1}$, both actions are equivalent.\\

We can prove that $(V,\oV,Z)$ is not an iLT in general, but we want to find constraints on the parameters such that this triplet becomes an iLT.   
As explained previously, to prove that some specializations of this triplet $(V,\oV,Z)$ are iLTs, it remains to show that 
there exist parameters $h_i$ such that relations \eqref{eq:Heung1} and \eqref{eq:Heung4} are satisfied. 

Acting with relations \eqref{eq:Heung1} and \eqref{eq:Heung4} on $\boldsymbol{z_i}$ and extracting the coefficients in front of $\boldsymbol{z_{i+1}}$, $\boldsymbol{z_{i}}$ and $\boldsymbol{z_{i-1}}$ lead to the constraints:
\begin{subequations}
    \begin{align}
 &   \zeta_{i}C^{(qR)}_{i+1}(\oalpha,\obeta,\ogamma,\odelta)\nu_{i+1}=(h_2+h_3\zeta_{i}+h_4\zeta_{i+1})C^{(qR)}_{i+1}(\alpha,\beta,\gamma,\delta)\,,\label{eq:zip}\\
 &   \zeta_{i} B^{(qR)}_{i}(\oalpha,\obeta,\ogamma,\odelta)=h_0+h_1\zeta_{i}+(h_2+(h_3+h_4)\zeta_{i})B^{(qR)}_{i}(\alpha,\beta,\gamma,\delta)\,,\label{eq:zi0}\\
  &   \zeta_{i} A^{(qR)}_{i-1}(\oalpha,\obeta,\ogamma,\odelta)\frac{1}{\nu_{i}}=(h_2+h_3\zeta_{i}+h_4\zeta_{i-1})A^{(qR)}_{i-1}(\alpha,\beta,\gamma,\delta)\,,\label{eq:zim}
\end{align}
\end{subequations}
and
\begin{subequations}
    \begin{align}
 &   \zeta_{i+1}C^{(qR)}_{i+1}(\alpha,\beta,\gamma,\delta)=(h_7+h_8\zeta_{i}+h_9\zeta_{i+1}){C}^{(qR)}_{i+1}(\oalpha,\obeta,\ogamma,\odelta)\nu_{i+1}\,,\label{eq:Zip}\\
 &   \zeta_{i} B^{(qR)}_{i}(\alpha,\beta,\gamma,\delta)=h_5+h_6\zeta_{i}+(h_7+(h_8+h_9)\zeta_{i}){B}^{(qR)}_{i}(\oalpha,\obeta,\ogamma,\odelta)\,,\label{eq:Zi0}\\
  &   \zeta_{i-1} A^{(qR)}_{i-1}(\alpha,\beta,\gamma,\delta)=(h_7+h_8\zeta_{i}+h_9\zeta_{i-1}){A}^{(qR)}_{i-1}(\oalpha,\obeta,\ogamma,\odelta)\frac{1}{\nu_{i}}\,.\label{eq:Zim}
\end{align}
\end{subequations}
We consider the previous relations as rational functions in terms of $q^i$.
Looking at the leading coefficients of \eqref{eq:zi0} and \eqref{eq:Zi0}, we deduce that 
\begin{align}
    h_1=h_6=0\,.
\end{align}
To simplify further computations, define the new parameters $s$, $h_{3,4}$,  $h_{4,3}$, $h_{8,9}$, $h_{9,8}$, $\mu_2$ and $\mu_7$ by
\begin{align}
 &  \sigma=-1/s-qs \alpha \beta\,,\quad h_{n,m}=h_n+qh_m \quad \text{for } (n,m)\in\{(3,4),(4,3),(8,9),(9,8)\} \,,\nonumber\\
  &  h_2=-h_{4,3}/\mu_2-h_{3,4}\alpha\beta\mu_2-\sigma(h_3+h_4)\,,\quad h_7=-h_{9,8}/\mu_7-h_{8,9}\alpha\beta\mu_7-\sigma(h_8+h_9)\,,\label{eq:parmu}
\end{align}
and the function $\psi$ by
\begin{align}
    \psi(\alpha,\beta,\gamma,\delta)=q(\alpha+\gamma  ) (1+\beta\delta)+q(\beta+\gamma  ) (\alpha +\delta )\,.  
\end{align}
Relation \eqref{eq:zip} allows for finding an expression of $\nu_{i}$ which can be used in \eqref{eq:Zip} to obtain the following relation:
\begin{align}
     \zeta_{i-1} \zeta_{i}
     =(h_2+h_3\zeta_{i-1}+h_4\zeta_{i})(h_7+h_8\zeta_{i-1}+h_9\zeta_{i})\,,\
\end{align}
which is given explicitly by 
\begin{align}
   & \frac{1}{s^2}(q^{i}-sq)(q^{i}-s)(\alpha\beta s q^{i}-1)(\alpha\beta s q^{i+1}-1)\label{eq:Cons1}\\
   =&\frac{1}{\mu_2\mu_7}(q^{i}-\mu_2)(q^{i}-\mu_7)(\alpha\beta\mu_2 h_{3,4}q^{i}-h_{4,3})(\alpha\beta\mu_7 h_{8,9}q^{i}-h_{9,8})\,.\nonumber
\end{align}
Computing the leading coefficient, we obtain that
\begin{align}\label{eq:h43}
 h_{8,9}=\frac{q}{h_{3,4}}\,.
\end{align}
Similar computations from \eqref{eq:zim} and \eqref{eq:Zim} yield
\begin{align}
  h_{4,3}=\frac{\odelta q}{\delta h_{3,4}}\,,\qquad  h_{9,8}=\frac{\delta h_{3,4}}{\odelta}
    \,. 
\end{align}
Using all the previous relations, 
the leading coefficients of \eqref{eq:zi0} and \eqref{eq:Zi0} provide the following expressions for $h_0$ and $h_5$:
\begin{subequations}\label{eq:h0h5}
\begin{align}
    &h_0=\psi(\oalpha,\obeta,\ogamma,\odelta)  -\frac{\delta h^2_{3,4}+q \odelta}{\delta h_{3,4}(1+q)}\psi(\alpha,\beta,\gamma,\delta)\,,\\
    &h_5=\psi(\alpha,\beta,\gamma,\delta)  -\frac{\delta h^2_{3,4} +q \odelta}{\odelta h_{3,4}(1+q)}\psi(\oalpha,\obeta,\ogamma,\odelta)\,.
\end{align}
\end{subequations}
Understanding relation \eqref{eq:Cons1} as an equality between two polynomials in $q^i$, their roots must be equal which leads to the equality between the sets:
\begin{align}
    \left\{ sq,s,\frac{1}{\alpha\beta s},\frac{1}{\alpha\beta s q}\right\}= \left\{ \mu_2,\mu_7, \frac{\odelta q}{ h_{3,4}^2 \alpha\beta\delta \mu_2}, \frac{\delta h_{3,4}^2}{\alpha\beta\odelta \mu_7 q} \right\}\,.\label{eq:2sets}
\end{align}
The product of the 4 elements of the set on the left is equal to the product of the 4 elements of the set on the right. Therefore, comparing these two sets leads only to three independent equations. We can determine for example $\mu_2$, $\mu_7$ and $\odelta$. We obtain 4!=24 different solutions.

In the following, we suppose that the parameters $\alpha$, $\beta$ and $\delta$ are generic and we look for the parameters $\oalpha$, $\odelta$, $\sigma$ and the function $\nu_i$ such that the previous relations hold. With mathematical software, we test the 24 different solutions in the remaining equations and solve in terms of $h_{3,4}$, $\widetilde{\alpha}$ and $s$. 
We consider only the solutions which provide different expressions for $\oV$ (by different, we mean $\oV$ that are not proportional and that are not related by transformations on the remaining free parameters). 
In the end, we find different possibilities for the uplet $(\oalpha,\odelta,\mu_2,\mu_7,h_{3,4},s)$.
However, asking that the conditions \eqref{eq:cond1} and \eqref{eq:cond2} be satisfied, we arrive to only one solution:
\begin{align}\label{eq:parI}
  (\oalpha,\odelta,\mu_2,\mu_7,h_{3,4},s)=  (\alpha,\alpha/(\beta\delta),s,1/(\alpha\beta s),1/(\beta\delta s),s)\,.
\end{align}
We detail in the following subsection, the iLT obtained from this solution, and show that the associated overlap coefficients are proportional to the Wilson rational functions.

\subsection{Wilson rational functions and irreducible Leonard trio  \label{sec:LTwilson}}

From now on, one considers that relation \eqref{eq:parI} holds. 
We will show that the overlap coefficients associated to this LT are the Wilson rational functions.

\paragraph{Basis $\boldsymbol{z_i}$.}
The action of $Z$ on $\boldsymbol{z_i}$ is diagonal and is given by:
\begin{align}
    Z\boldsymbol{z_i}&=(\lambda^{(qR)}(i;\alpha\beta)+\sigma)\boldsymbol{z_i}=(q^{-i}-1/s)(1-\alpha\beta s q^{i+1})
\boldsymbol{z_i}\,,
  \end{align}
  where $\lambda^{(qR)}(i;\alpha\beta)$ is defined in Appendix \ref{sec:qRacah}.
  The normalization $\nu_i$ becomes
  \begin{equation}
      \nu_i=\frac{\beta(s-q^{i})(\delta-\alpha q^{i})}{(1-\beta\delta q^{i})(1-\alpha\beta s q^{i})}\,,
  \end{equation}
and the tridiagonal actions of $V$ and $\oV$ are
 \begin{align}
    V\boldsymbol{z_i}&=C^{(qR)}_{i+1}(\alpha,\beta,\gamma,\delta)\boldsymbol{z_{i+1}}+B^{(qR)}_i(\alpha,\beta,\gamma,\delta)\boldsymbol{z_i}
   +A^{(qR)}_{i-1}(\alpha,\beta,\gamma,\delta)\boldsymbol{z_{i-1}}\,,\label{eq:Vz1}\\
    \oV\boldsymbol{z_i}&=C^{(qR)}_{i+1}(\alpha,\beta,\gamma,\alpha/(\beta\delta))\frac{\beta(s-q^{i+1})(\delta-\alpha q^{i+1})}{(1-\beta\delta q^{i+1})(1-\alpha\beta s q^{i+1})}\boldsymbol{z_{i+1}}+B^{(qR)}_i(\alpha,\beta,\gamma,\alpha/(\beta\delta))\boldsymbol{z_i}\nonumber\\
    &+A^{(qR)}_{i-1}(\alpha,\beta,\gamma,\alpha/(\beta\delta))\frac{(1-\beta\delta q^{i})(1-\alpha\beta s q^{i})}{\beta(s-q^{i})(\delta-\alpha q^{i})}\boldsymbol{z_{i-1}} \,,\label{eq:Vtz1}
\end{align}
where we recall that the functions $A^{(qR)}_i$, $B^{(qR)}_i$ and $C^{(qR)}_i$ are defined in Appendix \ref{sec:qRacah}.

\paragraph{Algebraic Heun operators.}
In this case, the parameters entering in the definition of the algebraic Heun operators can be given explicitly, using \eqref{eq:parmu} and \eqref{eq:h0h5}, and relations \eqref{eq:Heung1} and \eqref{eq:Heung4} become 
\begin{align}
    \oV Z&=h_0 \un+\frac{(1-\alpha\beta s^2 )(1-\alpha\beta q^2s^2)}{\beta\delta s^2(1+q )} V +\frac{1}{1-q^2}\left( \frac{\alpha q s}{\delta}[Z,V]_q+\frac{1}{\beta \delta s}[V,Z]_q\right)  \,,\label{eq:Heunovz}\\
    Z V&=h_5 \un+\frac{\delta(1-\alpha\beta s^2)(1-\alpha\beta q^2s^2)}{\alpha s^2(1+q )}\oV +\frac{1}{1-q^2}\left(  \frac{\delta}{s\alpha}[Z,\oV]_q+q\beta\delta s[\oV,Z]_q\right)  \,,  \label{eq:Heunzv} 
\end{align}
with $[Z,V]_q=ZV-qVZ$ and 
\begin{align}
  &   h_0=\psi(\alpha,\beta,\gamma,\alpha/(\beta\delta)) 
  - \frac{1+q s^2\alpha\beta}{\beta\delta s(1+q)}\psi(\alpha,\beta,\gamma,\delta)\,,\label{eq:h0W}\\
 & h_5=\psi(\alpha,\beta,\gamma,\delta) 
 - \frac{\delta(1+q s^2\alpha\beta)}{\alpha s(1+q)}\psi(\alpha,\beta,\gamma,\alpha/(\beta\delta))\,.
 \label{eq:h5W}
\end{align}
In addition to the two relations \eqref{eq:Heunovz} and \eqref{eq:Heunzv}, there are two other similar relations given in the following proposition. 
\begin{prop} The following relations hold
    \begin{align}
    (Z+\rho\un) \oV&=\widetilde{h}_0 \un+\frac{q(1-s^2\alpha\beta)^2}{\beta\delta s^2 (1+q)} V+\frac{1}{1-q^2}\left(  \frac{s\alpha}{\delta}[Z,V]_q+\frac{q}{\beta\delta s}[V,Z]_q\right)   \,,\label{eq:Heunzov2}\\
     V(Z+\rho\un) &=\widetilde{h}_5 \un+\frac{q\delta(1-s^2\alpha\beta)^2}{\alpha s^2(1+q)} \oV +\frac{1}{1-q^2}\left(\frac{q\delta}{s\alpha}[Z,\oV]_q+ s\beta\delta [\oV,Z]_q\right)  \,,\label{eq:Heunvz2}
\end{align}
where $\rho=(1-q )(1/s-\alpha\beta s )$ and
\begin{align}
 &   \widetilde{h}_0=\psi(\alpha,\beta,\gamma,\alpha/(\beta\delta)) 
  - \frac{q+ s^2\alpha\beta}{\beta\delta s(1+q)}\psi(\alpha,\beta,\gamma,\delta)\,,\\
 & \widetilde{h}_5=\psi(\alpha,\beta,\gamma,\delta) 
 - \frac{\delta(q+ s^2\alpha\beta)}{\alpha s(1+q)}\psi(\alpha,\beta,\gamma,\alpha/(\beta\delta))\,.
\end{align}
\end{prop}
\proof Relations of this proposition follow from \eqref{eq:Heunovz} and \eqref{eq:Heunzv}, see Corollary \ref{cor:otherHeun}. \endproof
The statement of the previous proposition implies that $(Z+\rho\un) \oV$ is an algebraic Heun operator associated to the LP $(V,Z)$, and $V(Z+\rho\un)$ is an algebraic Heun operator associated to the LP $(\oV,Z)$.
It is one way to obtain the bispectral properties of the biorthogonal partner of $w_n(x)$ (see the comment at the end of this subsection).

As stated in Proposition \ref{prop:GEVPg}, the bispectral relations for $w_n(x)$ are obtained by the actions of $Z$, $\oV Z$ on $\boldsymbol{v_n}$ and the actions of $Z$, $Z V $ on $\boldsymbol{\ov_x}$. The next two paragraphs provide these actions.

\paragraph{Basis $\boldsymbol{v_n}$.} 
As explained in Subsection \ref{sec:LPbp}, the vectors $\boldsymbol{v_n}$ diagonalizing $V$, can be expressed in terms of $\boldsymbol{z_i}$ by, for $n=0,1,\dots,N$,
\begin{align}\label{eq:vnzi}
    \boldsymbol{v_n}=\sum_{i=0}^N R^{(qR)}_i(n;\alpha,\beta,\gamma,\delta)\ \boldsymbol{z_{i}}\,,
\end{align}
where $R^{(qR)}_i$ are the $q$-Racah polynomials defined in Appendix \ref{sec:qRacah}. Using the recurrence relation  of the $q$-Racah polynomials \eqref{eq:recuR}, one can show that
\begin{align}
    V \boldsymbol{v_n}=\lambda^{(qR)}(n;\gamma\delta)\boldsymbol{v_n}\,.
\end{align}
The difference relation of the $q$-Racah polynomials \eqref{eq:diffR} implies that
\begin{align}
  &   Z \boldsymbol{v_n} =A^{(qR)}_{n}(\gamma,\delta,\alpha,\beta) \boldsymbol{v_{n+1}}+\big(B^{(qR)}_{n}(\gamma,\delta,\alpha,\beta)+\sigma\big)\boldsymbol{v_n}+C^{(qR)}_{n}(\gamma,\delta,\alpha,\beta)\boldsymbol{v_{n-1}}\,,
\end{align}
where we recall that $\sigma=-\alpha\beta s q-1/s$.
The Heun operators associated to the LP $(V,Z)$ are tridiagonal in the basis $\boldsymbol{v_n}$.
Therefore, from the expression \eqref{eq:Heunovz}, the action of $\oV Z $ on $\boldsymbol{v_n}$ can be computed:
\begin{align}
    \oV Z \boldsymbol{v_n}&=\frac{1}{q \beta \delta s}\lambda^{(qR)}(n;\alpha\beta\gamma\delta s^2 q^2)A^{(qR)}_n(\gamma,\delta,\alpha,\beta)\boldsymbol{v_{n+1}}
    +\frac{\alpha s q}{  \delta}
\lambda^{(qR)}(n;\frac{\gamma\delta}{\alpha  \beta  q^2 s^2})C^{(qR)}_n(\gamma,\delta,\alpha,\beta)\boldsymbol{v_{n-1}}\nonumber\\
    &+\left(h_0-\frac{\alpha}{\delta} (1+q) \lambda^{(qR)}(n;\gamma\delta)-\frac{\sigma \lambda^{(qR)}(n;\gamma\delta)}{\beta  \delta  \left(1+q \right)}B^{(qR)}_n(\gamma,\delta,\alpha,\beta)\right)\boldsymbol{v_{n}}\,,
    \end{align}
    where $h_0$ is given by \eqref{eq:h0W}.

\paragraph{Basis $\boldsymbol{\ov_x}$.}
Similarly, the vectors $\boldsymbol{\ov_x}$, diagonalizing $\oV$ are given by, for $x=0,1,\dots,N$,
\begin{align}\label{eq:ovxvi7}
    \boldsymbol{\ov_x}=\sum_{i=0}^N \frac{(\beta \delta s)^i(q/s,q\alpha/\delta;q)_i}{(\beta \delta q,\alpha\beta s q;q)_i} R^{(qR)}_i(x;\alpha,\beta,\gamma,\alpha/(\beta\delta))\boldsymbol{z_{i}}\,,
\end{align}
where again $R^{(qR)}_i$ are the $q$-Racah polynomials defined in Appendix \ref{sec:qRacah}. 
The normalization appearing in the front of the $q$-Racah polynomials in the previous relation is explained by the factor $\nu_i$ in relation \eqref{eq:Vtz1}. This factor is $\prod_{k=1}^{i}\nu_k$.

As previously, using the properties of the $q$-Racah polynomials, one can show that
\begin{align}
    \oV \boldsymbol{\ov_x}=\lambda^{(qR)}(x;\alpha\gamma/(\beta\delta))\boldsymbol{\ov_x}\,,
\end{align}
and that
\begin{align}
  &   Z \boldsymbol{\ov_x} =A^{(qR)}_x(\gamma,\alpha/(\beta\delta),\alpha,\beta) \boldsymbol{\ov_{x+1}}+(B^{(qR)}_x(\gamma,\alpha/(\beta\delta),\alpha,\beta) +\sigma)\boldsymbol{\ov_{x}}+C^{(qR)}_x(\gamma,\alpha/(\beta\delta),\alpha,\beta) \boldsymbol{\ov_{x-1}}\,.
\end{align}
Using these results and relation \eqref{eq:Heunzv}, one proves that the action of $ZV$ is given by
\begin{align}
    Z V \boldsymbol{\ov_x}&=\beta \delta s \lambda^{(qR)}(x;\frac{\gamma}{\beta^2\delta s^2})A^{(qR)}_x(\gamma,\alpha/(\beta\delta),\alpha,\beta)\boldsymbol{\ov_{x+1}}+\frac{\delta}{\alpha s}\lambda^{(qR)}(x;\frac{\alpha^2\gamma s^2}{\delta})C^{(qR)}_x(\gamma,\alpha/(\beta\delta),\alpha,\beta)\boldsymbol{\ov_{x-1}}\nonumber\\
    &+\left(h_5-(1+q)\beta\delta\lambda^{(qR)}(x;\alpha\gamma/(\beta\delta))
    -\frac{\sigma \delta\lambda^{(qR)}(x;\alpha\gamma/(\beta\delta))}{ \alpha (1+q)}
B^{(qR)}_x(\gamma,\alpha/(\beta\delta),\alpha,\beta)\right)\boldsymbol{\ov_{x}}\,.
\end{align}

\paragraph{Dual eigenvectors.} We choose the scalar product $\langle\cdot,\cdot\rangle$ such that, for $i,j=0,1,\dots,N$,
\begin{align}
    \langle \boldsymbol{z_i},\boldsymbol{z_j}\rangle =\delta_{i,j}\,.
\end{align}
It is always possible to choose such a scalar product since the vector space $\V$ is finite dimensional.
For this choice, the dual vectors $\boldsymbol{\ov^\star_x}$ of $\boldsymbol{\ov_x}$ are given explicitly by
\begin{align}\label{eq:vtszi}
\boldsymbol{\ov^\star_x}  =
 M_{\gamma,\alpha/(\beta\delta),\alpha,\beta} \Omega_{x;\gamma,\alpha/(\beta\delta),\alpha,\beta}  \sum_{i=0}^N  \frac{(\beta \delta q,\alpha\beta s q;q)_i}{(\beta \delta s)^i(q/s,q\alpha/\delta;q)_i} \Omega_{i;\alpha,\beta,\gamma,\alpha/(\beta\delta)} R^{(qR)}_i(x;\alpha,\beta,\gamma,\alpha/(\beta\delta))\boldsymbol{z_i}\,,
\end{align}
where $M$ and $\Omega$ are defined by \eqref{eq:qRweight}. 
We choose the parameters such that the weight is positive. 
To prove that $\langle\boldsymbol{\ov^\star_x},\boldsymbol{\ov_y}\rangle=\delta_{x,y}$, we use the orthogonality relation satisfied by the $q$-Racah polynomials \eqref{eq:orthoqR}.

\paragraph{The Wilson rational functions.}
Using relations \eqref{eq:vnzi} and \eqref{eq:vtszi}, the following expression of the overlap coefficients $w_{n}(x)$ in terms of the $q$-Racah polynomials is obtained:
\begin{align}\label{eq:wWilson}
w_{n}(x)=\langle \boldsymbol{\ov_x^\star},\boldsymbol{v_n}\rangle& =M_{\gamma,\alpha/(\beta\delta),\alpha,\beta} \Omega_{x;\gamma,\alpha/(\beta\delta),\alpha,\beta}\\
&\times \sum_{i=0}^N \Omega_{i;\alpha,\beta,\gamma,\alpha s}R^{(qR)}_i(x;\alpha,\beta,\gamma,\alpha/(\beta\delta))R^{(qR)}_i(n;\alpha,\beta,\gamma,\delta)\,,\nonumber
\end{align}
where the functions $M$ and $\Omega$ are defined by \eqref{eq:qRweight}. These functions are closely related to the Wilson rational functions $W_n$, whose definition is recalled in Appendix \ref{sec:Wilson}, as stated precisely in the following proposition.
\begin{prop} \label{pro:wnxWilson}
    The overlap coefficients $w_n(x)$ defined by \eqref{eq:wWilson} can be expressed as follows:
     \begin{align}
       w_n(x)=&M_{\alpha,1/(q\alpha s),\gamma,\alpha/(\beta\delta)}\ \Omega_{x;\alpha,\gamma/(\beta\delta),\gamma,1/(s q \gamma)}\nonumber\\
     & \times 
      \frac{(\delta/(\alpha s),\gamma q/\beta;q)_n(\alpha\beta s q)^n } {(\beta\delta q,\alpha\gamma s q^2;q)_n}W_n(x;\alpha\gamma s q,\gamma\delta,\alpha\gamma/(\beta\delta),1/\gamma,1/\alpha,\beta/\gamma)\,,\label{eq:pr1}   
    \end{align}
    where $W_n$ is the Wilson rational functions defined in terms of the $q$-hypergeometric function  ${}_{10}\phi_9$ given in \eqref{eq:qRacahW}.
\end{prop}
\proof 
Using the results of Proposition \ref{prop:GEVPg}  with the different representations obtained in the previous paragraphs, one can compute the recurrence relation and the difference equation satisfied by $w_n(x)$.
The comparison with equations \eqref{eq:GEVPW1} satisfied by the normalized Wilson rational functions (r.h.s. of \eqref{eq:pr1}) proves that these normalized Wilson rational functions are equal to $w_n$ up to a proportionality constant independent of $n$ and $x$. It remains to prove that this constant is 1. To do that, let us compute $w_0(0)$, starting from its expression \eqref{eq:wWilson}:
\begin{align}
 w_0(0)&=\frac{ M_{\gamma,\alpha/(\beta \delta),\alpha,\beta}}{1-\alpha\beta q} \sum_{i=0}^N 
 \frac{(\alpha q,\gamma q,\alpha\beta q,\alpha\beta s q;q)_i(1-\alpha\beta q^{2i+1})}{(\alpha\gamma qs)^i(q,\alpha\beta q/\gamma ,\beta q,q/s;q)_i}=\frac{(q\beta\delta,1/(\alpha s);q)_N}{(\beta\delta/\alpha,q/s;q)_N}\,.\label{eq:pr2}
\end{align}
The last equality is obtained from the orthogonality relation of the $q$-Racah polynomials \eqref{eq:orthoqR} for $n=m=0$ and for $\alpha\to \gamma,\ \beta\to\alpha s,\ \gamma \to \alpha,\ \delta\to\beta$.
Therefore, the expressions \eqref{eq:wWilson} and the r.h.s. of \eqref{eq:pr1} for $n=x=0$ coincide which concludes the proof.
\endproof
The result of the previous proposition provides an expression of the Wilson rational functions in terms of the $q$-Racah polynomials. Indeed, using the explicit expressions of $M$ and $\Omega$, one can show that relation \eqref{eq:pr1} is equivalent to 
\begin{align}
&W_n(x;\alpha\gamma s q,\gamma\delta,\alpha\gamma/(\beta\delta),1/\gamma,1/\alpha,\beta/\gamma) =\frac{(\alpha q/\delta,\alpha\gamma s q^2;q)_x } {(1/(\beta\delta s),\gamma q/\beta;q)_x(\alpha\beta s q)^x}
      \frac{(\beta\delta q,\alpha\gamma s q^2;q)_n } {(\delta/(\alpha s),\gamma q/\beta;q)_n(\alpha\beta s q)^n }\nonumber\\
        &\qquad\qquad \times M_{\gamma,\alpha s,\alpha,\beta}\sum_{i=0}^N \Omega_{i;\alpha,\beta,\gamma,\alpha s} R^{(qR)}_i(x;\alpha,\beta,\gamma,\alpha/(\beta\delta))R^{(qR)}_i(n;\alpha,\beta,\gamma,\delta)\,.\label{eq:WRR}
\end{align}

\paragraph{Biorthogonal partner.}
Using the orthogonality of the $q$-Racah polynomials \eqref{eq:orthoqR}, we can show that the dual vectors of $\boldsymbol{v_n}$ are given by:
\begin{align}
    \boldsymbol{v^\star_n}=M_{\gamma,\delta,\alpha,\beta}\Omega_{n;\gamma,\delta,\alpha,\beta}\sum_{i=0}^N\Omega_{i;\alpha,\beta,\gamma,\delta}
    R^{(qR)}_i(n;\alpha,\beta,\gamma,\delta)\   \boldsymbol{z_i}\,,
\end{align}
which allows us to obtain the following expression for the biorthogonal partner (using expression \eqref{eq:ovxvi7} for $\boldsymbol{\ov_x}$):
\begin{align}
   \ow_n(x)= \langle \boldsymbol{v^\star_n},\boldsymbol{\ov_x}\rangle&= M_{\gamma,\delta,\alpha,\beta}\Omega_{n;\gamma,\delta,\alpha,\beta}\\
 &\times\sum_{i=0}^N \Omega_{i;\alpha,\beta,\gamma,1/(\beta s)} R^{(qR)}_i(x;\alpha,\beta,\gamma,\alpha/(\beta\delta))R^{(qR)}_i(n;\alpha,\beta,\gamma,\delta)\,.\nonumber
\end{align}
Comparing with the expression \eqref{eq:wWilson}, we see that 
\begin{align}
  \ow_n(x)  =M_{\gamma,\delta,\alpha,\beta\delta/(\alpha q)}\frac{\Omega_{n;\gamma,\delta,\alpha,\beta}}{ \Omega_{x;\gamma,\alpha/(\beta\delta),\alpha,\beta}}\  w_n(x)\big|_{s\to 1/(\alpha\beta s)}\,.
\end{align}
We deduce that the biorthogonal partners are also bispectral. The explicit expressions of these relations are not given here but can be computed easily from the ones of $w_n(x)$.
Let us emphasize that these bispectral relations could also be computed using computations similar to the ones done in the proof of Proposition \ref{prop:GEVPg}
and using the results of Proposition \ref{pro:bpw}. 

\paragraph{Generalized eigenvalue problem.}
The GEVP for the Wilson rational function are given in Appendix \ref{sec:Wilson}, we repeat them here  in the natural notations that emerge from our construction. The recurrence GEVP for the 
 Wilson rational functions $W_n(x)=W_n(x;\alpha\gamma s q,\allowbreak \gamma\delta,\allowbreak 
\alpha\gamma/(\beta\delta),\allowbreak 1/\gamma,\allowbreak 1/\alpha,\allowbreak 
\beta/\gamma)$ is: 
\begin{align}\label{eq:GEVPWn}
 &\X_{n+1,n}\, W_{n+1}(x)  +\X_{n,n}\, W_{n}(x)+ \X_{n-1,n} W_{n-1}(x)\nonumber\\
=&(1-q^x)(q^{-x}-\alpha\gamma q/(\beta\delta))\Big( \Z_{n+1,n}\, W_{n+1}(x)+ \Z_{n,n}\, W_{n}(x)    + \Z_{n-1,n}\, W_{n-1}(x)\Big)\,,
\end{align}
where 
\begin{align}
& \Z_{n+1,n}=A^{(qR)}_n(\gamma,\delta,\alpha,\beta)\frac{q(\alpha s-\delta q^{n})(\beta-\gamma q^{n+1})}{(1-\alpha\gamma s q^{n+2})(1-\beta\delta q^{n+1})}\,,\\
 &\Z_{n-1,n}=C^{(qR)}_n(\gamma,\delta,\alpha,\beta)\frac{(1-\alpha\gamma s q^{n+1})(1-\beta\delta q^{n})}{q(\alpha s-\delta q^{n-1})(\beta-\gamma q^{n})}\,,\\
&\Z_{n,n}=B^{(qR)}_n(\gamma,\delta,\alpha,\beta) -\alpha\beta s q-1/s\,,
\end{align}
and
\begin{align}
& \X_{n+1,n}=A^{(qR)}_n(\gamma,\delta,\alpha,\beta)\frac{(\alpha s- \delta q^{n})(\beta-\gamma q^{n+1})(1-\beta\delta sq^{n+1})}{\beta\delta s q^{n}(1-\beta\delta q^{n+1})}
\,,\\
& \X_{n-1,n}=C^{(qR)}_n(\gamma,\delta,\alpha,\beta)\frac{(1-\alpha\gamma sq^{n+1})(1-\beta\delta q^{n})(\beta s -\gamma q^{n})}{\beta\delta s q^{n}(\beta-\gamma q^{n})}\,,\\
& \X_{n,n}= -\X_{n+1,n}-\X_{n-1,n}\,.
\end{align}
The difference GEVP for the Wilson rational function is obtained from its recurrence relation \eqref{eq:GEVPWn}, using the duality property satisfied by $W_n$ given by
\begin{align}
 W_n(x;\alpha\gamma s q,\allowbreak \gamma\delta,\allowbreak 
\alpha\gamma/(\beta\delta),\allowbreak 1/\gamma,\allowbreak 1/\alpha,\allowbreak 
\beta/\gamma)=
W_x(n;\alpha\gamma s q,\alpha\gamma/(\beta\delta), \gamma\delta
,\allowbreak 1/\gamma,\allowbreak 1/\alpha,\allowbreak 
\beta/\gamma)\,.
\end{align}
This duality consists in exchanging the variable $x$ and the degree $n$ with the change of parameter $\delta\to \alpha/(\beta\delta)$.

%%%%%%%%%%%%%%%%%%%%%%%%%%%%%%%%%%%%%
%%%%%%%%%%%%%%%%%%%%%%%%%%%%%%%%%%%%%
%%%%%%%%%%%%%%%%%%%%%%%%%%%%%%%%%%%%%
\section{Limits of the Wilson rational functions \label{sec:wrs}}

It is well-established that  other rational functions of interest and, in particular, the polynomials of the $q$-Askey scheme can be obtained from the Wilson rational function (see \textit{e.g.} \cite{GM,CTVZ}).

\subsection{$q$-Racah polynomials \label{ssec:poly}}

The $q$-Racah polynomials, which are at the top of the $q$-Askey scheme,  are obtained by considering the limit $s\to 0$, and the other polynomials are also reached using known limits (see \textit{e.g.} \cite{Koekoek}).
This result for the $q$-Racah polynomials can be obtained from different perspectives. The limit can be performed directly on the ${}_{10}\phi_9$ to get the ${}_4\phi_3$ form \eqref{eq:qRacahRF} of the $q$-Racah polynomials.
A second way consists in taking the limit in relation \eqref{eq:GEVPWn} (after a global multiplication by $s$) to recover the usual recurrence relation \eqref{eq:recuR} of the $q$-Racah polynomials. Finally, we can also remark that, in this limit, the associated Leonard trio becomes $(V,\oV,-\un)$ where $(V,\oV)$ is a LP associated to the $q$-Racah polynomials. In the following, some details are given and, in particular, limits of relation \eqref{eq:WRR} are provided and reproduce some important formulas for these polynomials.  

 The limit $s\to 0$ of the ${}_{10}\phi_9$ form of the Wilson rational function $W_n(x;\alpha\gamma s q,\allowbreak \gamma\delta,\allowbreak 
\alpha\gamma/(\beta\delta),\allowbreak 1/\gamma,\allowbreak 1/\alpha,\allowbreak 
\beta/\gamma)$ gives $R^{(qR)}_n(x;\gamma/\beta,\beta\delta,\gamma,\alpha/(\beta\delta))$.
Using this result in 
\eqref{eq:WRR}, we get
\begin{align}
R^{(qR)}_n(x;\rr_1)&=\frac{(\delta q^{-x}/\alpha;q)_x}{(\gamma q/\beta;q)_x}\frac{(q^{-n}/(\beta\delta);q)_n}{(\gamma q/\beta;q)_n}\frac{(\beta q;q)_N(\alpha q)^N}{(\alpha\beta q^2;q)_N}\nonumber\\
&\times \sum_{i=0}^N\frac{(\alpha q,q^{-i}/\gamma,\alpha\beta q;q)_i}{(q,\alpha\beta q/\gamma,\beta q;q)_i}\frac{(1-\alpha\beta q^{2i+1})}{(1-\alpha \beta q)(\alpha q)^i}R^{(qR)}_i(x;\rr_2)R^{(qR)}_i(n;\rr_3)\,,\label{eq:Racahrel}
\end{align}
where $\rr_1=\gamma/\beta,\beta\delta,\gamma,\alpha/(\beta\delta)$, $\rr_2=\alpha,\beta,\gamma,\alpha/(\beta\delta),$ and $\rr_3=\alpha,\beta,\gamma,\delta$.  
The relation just proven between three $q$-Racah polynomials is called the Racah relation. This relation for $q\to 1$ appear in the study of the $6j$-symbol for the Lie algebra $\mathfrak{su}(2)$ (see \textit{e.g.} \cite{Mess}) and has been proven for generic parameters using the representation theory of the Racah algebra \cite{Icosi}. When the $6j$-symbol for the quantum group
$U_q(\mathfrak{su}(2))$ is instead considered, there exists a similar relation involving the $q$-Racah polynomials with
particular parameters \cite{KR} and their generalizations for generic parameters can be proven using the Askey--Wilson algebra \cite{CFGR}. The construction proposed here provides a new proof of this relation.

\subsection{$q$-Racah type rational functions}
A limit case of the Wilson rational functions and of the corresponding Leonard trio appear in the representation theory of $U_q(\mathfrak{su}(2))$, see \cite{GW2025}. Setting $\alpha=\beta$ and $s=\sigma/\alpha$ in the Wilson rational function \eqref{eq:WRR} and letting $\alpha \to 0$ we obtain a very-well-poised $_8\phi_7$-function, which can be transformed into a balanced $_4\phi_3$-function using Watson's transformation \cite[(III.18)]{GR},
\[
{\mathscr{R}}^{(3)}_n(x) = \frac{(\gamma \sigma q^2, q^{-n}/\delta;q)_n}{(\sigma q^{-n}/\delta, \gamma q;q)_n} {}_{4}\phi_3 \left({{ q^{-n},\;\gamma\delta q^{n+1},\;\sigma q,\;\sigma \delta}\atop
{\gamma \sigma q^{x+1} ,\;\sigma\delta  q^{-x},\; \delta q}}\;\Bigg\vert \; q;q\right)\,.
\]
In this limit the two $q$-Racah polynomials in \eqref{eq:WRR} both become dual $q$-Krawtchouk polynomials and \eqref{eq:WRR} reduces to a formula expressing ${\mathscr{R}}^{(3)}_n(x)$ as a sum of products of two dual $q$-Krawtchouk polynomials. This is directly related to results in \cite[Sect.2]{GW2025}, where a multiple of the $q$-Racah type rational function ${\mathscr{R}}^{(3)}_n(x)$ appears as the overlap coefficients between two bases of $U_q(\mathfrak{su}(2))$-modules consisting of dual $q$-Krawtchouk polynomials. The recurrence relations for the rational functions are obtained using the actions of two (almost) twisted primitive elements and a standard Cartan element, which, in the language of the current paper, form an irreducible Leonard trio:
\[
V=\pi_N(\widetilde X_{1,s}),\quad \widetilde V=\pi_N(\widetilde X_{v,t}),\quad Z=\pi_N(K^2),
\]
where we use notation from \cite{GW2025}.

\subsection{Reduced Leonard trio}

Some limits of the Wilson rational functions lead to some functions associated to a reduced Leonard trio. 
For example, by taking the limit $\alpha\to0$ in the Wilson rational functions \eqref{eq:WRR}, we obtain\footnote{There is a change of notation in comparison to \cite{CTVZ}, where this function is denoted by ${\overline{\mathscr{R}}}^{(1)}_n(x)$}
\begin{equation}\label{eq:R1}
{\mathscr{R}}^{(1)}_n(x)={}_{4}\phi_3 \left({{ q^{-n},\;\gamma\delta q^{n+1},\;q^{-x},\;\gamma s q}\atop
{\gamma q ,\;\gamma q/\beta,\; \beta \delta sq^{1-x}}}\;\Bigg\vert \; q;q\right)\,.
\end{equation}
Formula \eqref{eq:WRR} becomes in this case:
\begin{align}
\label{eq:sum-R1}
{\mathscr{R}}^{(1)}_n(x)=&\frac{(q^{-n}/(\beta\delta);q)_n}{(1/(\beta\delta s);q)_x(\gamma q/\beta;q)_n }\sum_{i=N-x}^N \frac{(q^{i+1}/s;q)_{N-i}(q^{-x},\beta\delta q^{N-x+1};q)_{i+x-N} }{(q;q)_{i+x-N}(\beta\delta s q^{-x})^{i+x-N}}  R^{(dqH)}_i(n;\rr)\,,
\end{align}
where  $R^{(dqH)}_i$ are the dual $q$-Hahn polynomials defined by \eqref{eq:dqH} and $\rr=\beta\delta,\gamma/\beta, q^{-N-1}.$
The corresponding GEVP is 
\begin{align}
 &\X_{n+1,n}\, {\mathscr{R}}^{(1)}_{n+1}(x)  +\X_{n,n}\, {\mathscr{R}}^{(1)}_{n}(x)+ \X_{n-1,n} {\mathscr{R}}^{(1)}_{n-1}(x)\nonumber\\
=& (q^{-x}-1)\Big( \Z_{n+1,n}\, {\mathscr{R}}^{(1)}_{n+1}(x)+ \Z_{n,n}\, {\mathscr{R}}^{(1)}_{n}(x)    + \Z_{n-1,n}\, {\mathscr{R}}^{(1)}_{n-1}(x)\Big)\,,
\end{align}
where 
\begin{align}
 \Z_{n+1,n}&=-\frac{\delta  q^{n +1} \left(\beta-\gamma  q^{n +1} \right) \left(1-\gamma  q^{n +1}\right) \left(1-\gamma  \delta  q^{n +1}\right)}{\left(1-\gamma  \delta  q^{2 n +1}\right) \left(1-\gamma  \delta  q^{2 n +2}\right)}
\,,\\
 \Z_{n-1,n}&=\frac{q\gamma\left(1-\delta  q^{n}\right) \left(1-q^{n}\right) \left(1-\delta  \beta  q^{n}\right)}{\left(1-\gamma  \delta  q^{2 n}\right) \left(1-\gamma  \delta  q^{2 n +1}\right)}
\,,\\
\Z_{n,n}&=\frac{s-1}{s}-\frac{\left(1-\gamma  q^{n +1}\right) \left(1-\gamma  \delta  q^{n +1}\right) \left(1-\beta  \delta  q^{n +1}\right)}{\left(1-\gamma  \delta  q^{2 n +1}\right) \left(1-\gamma  \delta  q^{2 n +2}\right)}+\frac{\gamma  \delta  q^{n+1} \left(1-q^{n}\right) \left(1-\delta  q^{n}\right)  \left(\beta -\gamma  q^{n}\right)}{\left(1-\gamma  \delta  q^{2 n}\right) \left(1-\gamma  \delta  q^{2 n +1}\right)},
\end{align}
and
\begin{align}
 \X_{n+1,n}&= -\frac{\left(\beta-\gamma  q^{n +1} \right) \left(1-\beta  \delta  s q^{n +1} \right) \left(1-\gamma  q^{n +1}\right) \left(1-\gamma  \delta  q^{n +1}\right)}{\beta  s \left(1-\gamma  \delta  q^{2 n +1}\right) \left(1-\gamma  \delta  q^{2 n +2}\right)}\,,\\
 \X_{n-1,n}&=\frac{q\gamma\left(1-\delta  \beta  q^{n}\right) \left(\gamma  q^{n}-\beta  s \right) \left(1-q^{n}\right) \left(1-\delta  q^{n}\right) }{\beta  s \left(1-\gamma  \delta  q^{2 n}\right) \left(1-\gamma  \delta  q^{2 n +1}\right)}\,,
\qquad \X_{n,n}= -\X_{n+1,n}-\X_{n-1,n}\,,
\end{align}
where $\rr=\beta\delta, \gamma/\beta,q^{-N-1}$. 

These functions satisfy a $R_I$ difference equation given by
\begin{align}
 &\widetilde{\X}_{x,x+1}\, \mathscr{R}^{(1)}_{n}(x+1)  +\widetilde{\X}_{x,x}\, \mathscr{R}^{(1)}_{n}(x)+ \widetilde{\X}_{x,x-1} \mathscr{R}^{(1)}_{n}(x-1)\nonumber\\
=& (1-q^{n})(q^{-n}-\gamma\delta q)\Big( \widetilde{Z}_{x,x}\, \mathscr{R}^{(1)}_{n}(x)    + \widetilde{\Z}_{x,x-1}\, \mathscr{R}^{(1)}_{n}(x-1)\Big)\,,
\end{align}
where 
\begin{align}
 &\widetilde{\Z}_{x,x-1}=\frac{\left(1-q^{x}\right) \left(\beta\delta-\gamma q^{x}\right)}{\beta\delta s- q^{x-1}}\,,\qquad
\widetilde{\Z}_{x,x}=\gamma  q^{x+1}-\frac{1}{s}\,,
\end{align}
and
\begin{align}
& \widetilde{\X}_{x,x+1}= \frac{\left(1-\gamma q^{x +1}\right) \left(\beta\delta s- q^{x}\right) \left(\beta-\gamma  q^{x +1} \right)}{  \beta s q^{x }}\,,\\
& \widetilde{\X}_{x,x-1}=\frac{ \left(1-q^{x}\right) \left(\beta\delta -\gamma  q^{x}\right) \left(\beta  s -\gamma  q^{x}\right)}{\beta s q^{x-1}}\,,\qquad
 \widetilde{\X}_{x,x}= -\widetilde{\X}_{x,x+1}-\widetilde{\X}_{x,x-1}\,.
\end{align}
 We can conclude that these functions are overlap coefficients for a rLT. We can show that the associated LP is of type dual $q$-Hahn. 

\subsection{Functions associated to LT}

Performing other limits, we can obtain some functions which are not associated to an iLT or a rLT.  There are the ones already discussed in Subsection \ref{ssec:poly} corresponding to the $q$-Racah polynomials. However, there exist other possibilities as discussed below.

By substituting $\beta\to \beta/s$ in the rational function  $\mathscr{R}^{(1)}_n$, and take the limit $s\to 0$. We recover the rational function studied in \cite{bernard2024meta}
\begin{align}
\mathscr{H}^{(1)}_n(x)={}_{3}\phi_2 \left({{ q^{-n},\;\gamma\delta q^{n+1},\;q^{-x}}\atop
{\gamma q ,\; \beta \delta q^{1-x}}}\;\Bigg\vert \; q;q\right)\,.
\end{align}
The corresponding recurrence GEVP is 
\begin{align}
 &\X_{n+1,n}\, {\mathscr{H}}^{(1)}_{n+1}(x)  +\X_{n,n}\, {\mathscr{H}}^{(1)}_{n}(x)+ \X_{n-1,n} {\mathscr{H}}^{(1)}_{n-1}(x)\nonumber\\
=& (q^{-x}-1)\Big( \Z_{n+1,n}\, {\mathscr{H}}^{(1)}_{n+1}(x)+ \Z_{n,n}\, {\mathscr{H}}^{(1)}_{n}(x)    + \Z_{n-1,n}\, {\mathscr{H}}^{(1)}_{n-1}(x)\Big)\,,
\end{align}
where 
\begin{align}
& \Z_{n+1,n}=-\frac{\left(1-\gamma  \delta q^{n +1}\right) \left(1-\gamma  q^{n +1}\right) \beta  \delta q^{n +1}}{\left(1-\gamma  \delta  q^{2 n +2}\right) \left(1-\gamma  \delta  q^{2 n +1}\right)}
\,,\\
 &\Z_{n-1,n}=-\frac{\left(1-\delta  q^{n}\right) \left(1-q^{n}\right)  \gamma  \delta  \beta q^{n+1}}{\left(1-\gamma  \delta  q^{2 n +1}\right) \left(1-\gamma  \delta  q^{2 n}\right)}\,,\qquad
 \Z_{n,n}=-\Z_{n+1,n}-\Z_{n-1,n}-1\,,
\end{align}
and
\begin{align}
& \X_{n+1,n}= -\frac{\left(1-\gamma  \delta  q^{n +1}\right) \left(1-\gamma  q^{n +1}\right) \left(1-\beta  \delta  q^{n +1}\right)}{\left(1-\gamma  \delta  q^{2 n +2}\right) \left(1-\gamma  \delta  q^{2 n +1}\right)}
\,,\\
& \X_{n-1,n}=\frac{  \left(1-\delta q^{n}\right) \left(1-q^{n}\right) \left(\beta-\gamma  q^{n} \right) \delta\gamma  q^{n+1}}{\left(1-\gamma  \delta  q^{2 n +1}\right) \left(1-\gamma  \delta q^{2 n}\right)}\,,\qquad
 \X_{n,n}= -\X_{n+1,n}-\X_{n-1,n}\,.
\end{align} 
These functions satisfy a difference relation of type $R_I$ given explicitly by
\begin{align}
 &\widetilde{\X}_{x,x+1}\, \mathscr{H}^{(1)}_{n}(x+1)  +\widetilde{\X}_{x,x}\, \mathscr{H}^{(1)}_{n}(x)+ \widetilde{\X}_{x,x-1} \mathscr{H}^{(1)}_{n}(x-1)\nonumber\\
=& (1-q^{n})(q^{-n}-\gamma\delta q)\Big( \widetilde{\Z}_{x,x}\, \mathscr{H}^{(1)}_{n}(x)    + \widetilde{\Z}_{x,x-1}\, \mathscr{H}^{(1)}_{n}(x-1)\Big)\,,
\end{align}
where 
\begin{align}
 &\widetilde{\Z}_{x,x-1}=\frac{\beta\delta(1-q^x)}{\beta\delta-q^{x-1}}\,,\qquad
\widetilde{\Z}_{x,x}=-1\,,\label{eq:eig1}\\
& \widetilde{\X}_{x,x+1}=(\beta\delta-q^x)(q^{-x}-\gamma q)\,,\\
& \widetilde{\X}_{x,x-1}=\delta q (\beta-\gamma q^x)(q^{-x}-1) \,,\qquad
 \widetilde{\X}_{x,x}= -\widetilde{\X}_{x,x+1}-\widetilde{\X}_{x,x-1}\,.
\end{align}
These functions are not associated to a rLT because we can see, from the representation of $Z$ given by  \eqref{eq:eig1}, that this matrix is not diagonalizable which contradicts the definition of a rLT.

\section{Outlook \label{sec:out}}

In this paper, we introduced the concept of Leonard trios and explained their fundamental connection with bispectral rational functions. This approach establishes an algebraic framework for the study of these functions, opening several promising avenues for future research. It offers a perspective that differs from the meta algebra program \cite{vinet2021unified, tsujimoto2024meta, bernard2024meta, CTVZ, CLMTVZ2026} while exhibiting common elements. A review integrating both points of view would certainly be warranted. It would also be of interest to build a bridge to connect the Leonard Trio approach presented here to the description of the Wilson rational functions based on rational Heun operators \cite{TVZ2023rational}.

The classification of the Leonard trio would be a way to obtain a complete scheme for the bispectral rational functions analogous to the ($q$-)Askey scheme for the classical bispectral polynomials.  As a first objective, the complete classification of the irreducible Leonard trios would be a significant milestone. In particular, utilizing two Leonard pairs of Bannai--Ito type should allow for the definition of an associated Bannai--Ito bispectral rational function which, to our knowledge, has not been yet studied.

To achieve a full classification, it should be essential to determine the algebraic relations between the operators $V$, $\oV$, and $Z$ that constitute a Leonard trio. This would generalize the well-known Askey--Wilson relations satisfied by both elements of a Leonard pair $(X,Y)$ \cite{Vidunas}. This Askey--Wilson algebra connects Leonard pairs to diverse physical and mathematical problems, such as the algebraic interpretation of $3nj$-symbols or th symmetry algebras for superintegrable models. We anticipate that the discovery of the algebra underlying Leonard trios will have similar striking applications.

The notion of tridiagonal pairs, introduced in \cite{ITT}, generalizes Leonard pairs by allowing for eigenspaces with dimensions greater than one. These pairs are closely associated with tridiagonal relations \cite{PT} and various special functions \cite{BVZ, CGT}. A similar generalization is possible for  Leonard trios, and we believe that such an extension will lead to the discovery of novel algebras and special functions.

Finally, the extension of the approach developed in this paper is paving the way to pursue studies of the multivariate bispectral rational functions initiated in \cite{GW2025}. Indeed, just as higher-rank Leonard pairs were defined and studied to explore multivariate bispectral polynomials \cite{IT,CZ}, similar generalizations should be possible for Leonard trios.

\appendix

%%%%%%%%%%%%%%%%%%%%%%%%%%%%%%%%
\section{$q$-Racah polynomials \label{sec:qRacah} } 
The finite families of orthogonal polynomials in the $q$-Askey scheme are defined in terms of $q$-hypergeometric functions \cite{GR}:
\begin{align}
{}_{r+1}\phi_r\left({{q^{-n},\;a_1,\; a_2,\; \cdots,\; a_{r}  }\atop
{b_1,\; b_2,\; \cdots,\; b_{r} }}\;\Bigg\vert \; q;z\right)=\sum_{k=0}^{n}
\frac{(q^{-n},a_1,a_2,\cdots,a_r;q)_k}{(b_1,b_2,\cdots,b_r,q;q)_k}z^k\,,
\end{align}
for $r,n$ non-negative integers,
and where the $q$-Pochhammer symbols are
\begin{align}
(b_1,b_2,\cdots,b_r;q)_k=(b_1;q)_k(b_2;q)_k\cdots (b_r;q)_k\,,\qquad (b_i;q)_k=(1-b_i)(1-qb_i)\cdots (1-q^{k-1}b_i)\,.
\end{align}

The $q$-Racah polynomials are defined in terms of the $q$-hypergeometric functions as follows, for $n$ a non-negative integer:
\begin{align}\label{eq:qRacahRF}
R^{(qR)}_n(x;\rr)={}_{4}\phi_3 \left({{ q^{-n},\; \alpha\beta q^{n+1},\; q^{-x},\;  \gamma\delta  q^{x+1} }\atop
{  \alpha q,\;\beta\delta q,\; \gamma q }}\;\Bigg\vert \; q;q\right)\,,
\end{align}
where $\rr=\alpha,\beta,\gamma,\delta$ and $\alpha=q^{-N-1}$, $\beta\delta=q^{-N-1}$ or $\gamma=q^{-N-1}$ with $N$ a non-negative integer. 

They satisfy the recurrence relation, for $n,x=0,1,\dots N$,
\begin{align}
\label{eq:recuR}
    &\lambda^{(qR)}(x;\gamma\delta)R^{(qR)}_n(x;\rr)=A^{(qR)}_{n}(\rr) R^{(qR)}_{n+1}(x;\rr)+B^{(qR)}_{n}(\rr) R^{(qR)}_{n}(x;\rr)+C^{(qR)}_{n}(\rr)R^{(qR)}_{n-1}(x;\rr)\,,
\end{align}
where  $\lambda^{(qR)}(x;a)=q^{-x}+aq^{x+1}\,, $
and
\begin{align}
&A^{(qR)}_{n}(\rr)=\frac{(1-\alpha q^{n+1})(1-\alpha\beta q^{n+1})(1-\beta\delta q^{n+1})(1-\gamma q^{n+1})}{(1-\alpha\beta q^{2n+1})(1-\alpha\beta q^{2n+2})}\,,\\
     &C^{(qR)}_{n}(\rr)=\frac{q(1-q^{n})(1-\beta q^{n})(\gamma-\alpha\beta q^{n})(\delta-\alpha q^{n})}{(1-\alpha\beta q^{2n})(1-\alpha\beta q^{2n+1})}\,,\\
     & B^{(qR)}_{n}(\rr)=-A^{(qR)}_{n}(\rr)-C^{(qR)}_{n}(\rr)+1+\gamma\delta q\,.\label{eq:BAC}
\end{align}
The $q$-Racah polynomials satisfy the duality relation
\begin{align}
    R^{(qR)}_n(x;\rr)=R^{(qR)}_x(n;\rr')\,,
\end{align}
where $\rr'=\gamma,\delta,\alpha,\beta$.
From this duality relation and the recurrence relation, we can deduce the difference equation satisfied by the $q$-Racah polynomials:
\begin{align}
\label{eq:diffR}
    &\lambda^{(qR)}(n;\alpha\beta)R^{(qR)}_n(x;\rr)=A_{x}(\rr') R_{n}^{(qR)}(x+1;\rr)+B_{x}(\rr') R_{n}^{(qR)}(x;\rr)+C_{x}(\rr') R_{n}^{(qR)}(x-1;\rr)\,.
\end{align}
The $q$-Racah polynomials are orthogonal:
\begin{align}\label{eq:orthoqR}
 \sum_{x=0}^N\Omega_{x;\rr'} R^{(qR)}_{m}(x;\rr)R^{(qR)}_n(x;\rr)=\frac{1}{M_\rr
\Omega_{n;\rr} }\delta_{m,n}\,,
\end{align}
where 
\begin{align}\label{eq:qRweight}
  & \Omega_{n;\rr}=\frac{(\alpha q,\gamma q, \beta\delta q,\alpha \beta q; q)_n}{(q,\alpha\beta q/\gamma,\alpha q/\delta,\beta q;q)_n}\ \frac{1-\alpha\beta q^{2n+1}}{(\gamma\delta q)^n(1-\alpha\beta q)} \,,\\
   &M_\rr=\begin{cases}
   \frac{(\gamma q /\beta,\delta q;q)_N}{(1/\beta,\gamma\delta q^2;q)_N}&\text{if } \alpha=q^{-N-1}\\
    \frac{(\alpha\beta q/\gamma,\beta q;q)_N}{(\alpha\beta q^2,\beta/\gamma;q)_N}&\text{if }  \beta\delta=q^{-N-1}\\
       \frac{(\alpha q/\delta,\beta q;q)_N}{(\alpha\beta q^2,1/\delta;q)_N}&\text{if } \gamma=q^{-N-1}
   \end{cases}
\end{align}
The $q$-Racah polynomials are at the top of the $q$-Askey scheme for the finite families of classical orthogonal polynomials. The other families can be obtained by limits and specializations. In particular, one obtains the dual $q$-Hahn polynomials which are defined by, for $n$ a non-negative integer:
\begin{align}\label{eq:dqH}
R^{(dqH)}_n(x;\alpha,\beta,\gamma)={}_3\phi_2 \left({{q^{-n},\;q^{-x},\;\alpha\beta q^{x+1}}\atop
{\alpha q,\; \gamma q }}\;\Bigg\vert \; q,q\right)\,.
\end{align}
The following lemma provides another useful limit.
\begin{lemm} For $\gamma=q^{-N-1}$, these different relations hold:
\begin{align}\label{eq:lim0}
     R_n^{(qR)}(x;\alpha,\beta,\gamma,\delta/\beta)
   \mathop{\sim}_{\beta \to 0}\begin{cases}
    \displaystyle  \frac{(q^{-x};q)_n}{(\alpha q,\gamma q,\delta q;q)_n }\left(\frac{\gamma \delta q^{x+1}}{\beta}\right)^n\,,\quad\text{for } n\leq x\\
     \displaystyle  \frac{(q^{-n};q)_x}{(\alpha q,\gamma q,\delta q;q)_x }\left(\frac{\gamma \delta q^{x+1}}{\beta}\right)^x\,,\quad\text{otherwise}
   \end{cases}
\end{align}
\end{lemm}
\proof 
Using the expression as a sum of the $q$-hypergeometric expression of the $q$-Racah polynomials in the l.h.s. of \eqref{eq:lim0}, one gets
\begin{align}
   R_n^{(qR)}(x;\alpha,\beta,\gamma,\delta/\beta)&={}_4\phi_3 \left({{q^{-n},\; \alpha\beta q^{n+1}, \;q^{-x},\;\gamma\delta q^{x+1}/\beta}\atop
{\alpha q,\; \delta q ,\; \gamma q }}\;\Bigg\vert \; q,q\right)\\
   &\mathop{\sim}_{\beta \to 0} \; 
\sum_{k=0}^N \frac{(q^{-n},q^{-x};q)_k}{(q,\alpha q,\delta q,\gamma q;q)_k} 
\left(-\frac{q^{x+2}\gamma\delta }{\beta}\right)^k 
q^{\binom{k}{2}}\,.
\end{align}
Depending on the sign of $n-x$ the above sum stop at different stage and provide different results for the equivalent in $\beta$. The first result of the lemma is obtained after some manipulations of the $q$-Pochhammer symbols. 
\endproof

%%%%%%%%%%%%%%%%%%%%%%%%%%%%%%%%
%%%%%%%%%%%%%%%%%%%%%%%%%%%%%%%%
%%%%%%%%%%%%%%%%%%%%%%%%%%%%%%%%
\section{Wilson rational functions\label{sec:Wilson}}

The Wilson rational functions are defined in terms of the very-well-poised $q$-hypergeometric functions as follows, for $n$ a non-negative integer:
\begin{align}\label{eq:qRacahW}
W_n(x;a,b,c,d,e,f)={}_{10}\phi_9 \left({{a,\;q\sqrt{a},\;-q\sqrt{a},\; q^{-n},\; bq^{n+1},\; q^{-x},\;  c  q^{x+1} ,\; ad  ,\;ae,\;af }\atop
{\sqrt{a},\;- \sqrt{a},\; aq^{-n}/b ,\;aq^{n+1},\; a q^{-x}/ c  ,\;   aq^{x+1},\;  q/d,\; q/e,\; q/f }}\;\Bigg\vert \; q;q\right)\,,
\end{align}
with the relation between the parameters
\begin{align}
 b   c d e f  =1\,,
\end{align}
and we choose $d=q^{N+1}$, with $N$ a non-negative integer.

The Wilson rational functions $W_n(x)=W_n(x;a,b,c,d,e,f)$ satisfy the following recurrence GEVP \cite{Ros}, for $x,n$ non-negative integers: 
\begin{align}
 &\X_{n+1,n}\, W_{n+1}(x)  +\X_{n,n}\, W_{n}(x)+ \X_{n-1,n} W_{n-1}(x)\nonumber\\
=&(\lambda(x;c)-1-cq)\Big( \Z_{n+1,n}\, W_{n+1}(x)+ \Z_{n,n}\, W_{n}(x)    + \Z_{n-1,n}\, W_{n-1}(x)\Big)\,, \label{eq:GEVPW1}
\end{align}
where 
\begin{align}
& \Z_{n+1,n}=A^{(qR)}_n(1/d,bd,1/e,f/d)\frac{(a-bq^{n+1})(f-q^{n+1})}{(1-aq^{n+1})(1-bfq^{n+1})}\,,\\
 &\Z_{n-1,n}=C^{(qR)}_n(1/d,bd,1/e,f/d)\frac{(1-aq^{n})(1-bfq^{n})}{(a-bq^{n})(f-q^{n})}\,,\\
&\Z_{n,n}=B^{(qR)}_n(1/d,bd,1/e,f/d) -af-qbcf/a\,,
\end{align}
and
\begin{align}
& \X_{n+1,n}=A^{(qR)}_n(1/d,bd,1/e,f/d)\frac{(a-bq^{n+1})(f-q^{n+1})(c-aq^n)}{aq^n(1-bfq^{n+1})}
\,,\\
& \X_{n-1,n}=C^{(qR)}_n(1/d,bd,1/e,f/d)\frac{(1-aq^{n})(1-bfq^{n})(a-bcq^{n+1})}{abq^n(f-q^{n})}\,,\\
& \X_{n,n}= -\X_{n+1,n}-\X_{n-1,n}\,.
\end{align}
Using the symmetry relation
\begin{align}
    W_n(x)=W_x(n)|_{b\leftrightarrow c}\,,
\end{align}
we can obtain the difference GEVP satisfied by $W_n(x)$.

\paragraph{Acknowledgements:}
N.~Cramp\'e is partially supported by the international research project AAPT of the CNRS. L.~Vinet is funded in part by a Discovery Grant from the Natural Sciences and Engineering Research Council (NSERC) of Canada. Q. Labriet and L. Morey enjoy postdoctoral fellowships provided by this grant. W. Groenevelt and C. Wagenaar thank the Centre de Recherches Math\'ematiques at the Universit\'e de Montr\'eal for its hospitality.
%\printbibliography
%\bibliographystyle{utphys}
%\bibliography{bib_wilson}

\end{document}